\documentclass[reqno]{amsart}
\usepackage{graphicx}
\usepackage{pictex}
\usepackage{mathrsfs}
\usepackage{hyperref}
\begin{document}

%
%

\def\labelenumi{(\theenumi)}

\newtheorem{thm}{Theorem}[section]
\newtheorem{lem}[thm]{Lemma}
\newtheorem{conj}[thm]{Conjecture}
\newtheorem{cor}[thm]{Corollary}
\newtheorem{add}[thm]{Addendum}
\newtheorem{prop}[thm]{Proposition}
\theoremstyle{definition}
\newtheorem{defn}[thm]{Definition}
\theoremstyle{remark}
\newtheorem{rmk}[thm]{{\bf Remark}}
\newtheorem{example}[thm]{{\bf Example}}

\newcommand{\Cinfty}{{\mathbf C}_{\infty}}
\newcommand{\OmegaH}{\Omega/\langle H \rangle}
\newcommand{\hatOmegaHstar}{\hat \Omega/\langle H_{\ast}\rangle}
\newcommand{\SurfG}{\Sigma_g}
\newcommand{\TriangG}{T_g}
\newcommand{\TriangGOne}{T_{g,1}}
\newcommand{\ProjG}{\mathcal{P}_g}
\newcommand{\TeichG}{\mathcal{T}_g}
\newcommand{\CirclePackGTau}{\mathsf{CPS}_{g,\tau}}
\newcommand{\CrossRatio}{{\bf c}}
\newcommand{\CrossRatioGTau}{\mathcal{C}_{g,\tau}}
\newcommand{\CrossRatioOneTau}{\mathcal{C}_{1,\tau}}
\newcommand{\DeformGTau}{\mathcal{C}_{g,\tau}}
\newcommand{\Forget}{\mathit{forg}}
\newcommand{\Uniform}{\mathit{u}}
\newcommand{\Section}{\mathit{sect}}
\newcommand{\SLTwoC}{\mathrm{SL}(2,{\mathbf C})}
\newcommand{\SLTwoR}{\mathrm{SL}(2,{\mathbf R})}
\newcommand{\SUTwo}{\mathrm{SU}(2)}
\newcommand{\PSLTwoC}{\mathrm{PSL}(2,{\mathbf C})}
\newcommand{\GLTwoZ}{\mathrm{GL}(2,{\mathbf Z})}
\newcommand{\GLTwoC}{\mathrm{GL}(2,{\mathbf C})}
\newcommand{\PSLTwoR}{\mathrm{PSL}(2,{\mathbf R})}
\newcommand{\PGLTwoR}{\mathrm{PGL}(2,{\mathbf R})}
\newcommand{\PGLTwoZ}{\mathrm{PGL}(2,{\mathbf Z})}
\newcommand{\GLTwoR}{\mathrm{GL}(2,{\mathbf R})}
\newcommand{\PSLTwoZ}{\mathrm{PSL}(2,{\mathbf Z})}
\newcommand{\SLTwoZ}{\mathrm{SL}(2,{\mathbf Z})}
\newcommand{\nnn}{\noindent}
\newcommand{\MCG}{{\mathcal {MCG}}}
\newcommand{\MCSS}{{\mathcal S}^{\rm mc}_{\rm alg}}
\newcommand{\MMap}{{\bf \Phi}_{\mu}}
\newcommand{\HH}{{\mathbb H}^2}
\newcommand{\TT}{{\mathbb T}}
\newcommand{\X}{{\mathcal  X}}
\newcommand{\C}{{\mathscr C}}
\newcommand{\CC}{{\mathbf C}}
\newcommand{\RR}{{\mathbf R}}
\newcommand{\B}{{\mathcal  B}}
\newcommand{\G}{{\mathcal  G}}
\newcommand{\R}{{\mathcal  R}}
\newcommand{\Q}{{\mathbf Q}}
\newcommand{\ZZ}{{\mathbf Z}}
\newcommand{\PL}{{\mathscr {PL}}}
\newcommand{\GP}{{\mathcal {GP}}}
\newcommand{\GT}{{\mathcal {GT}}}
\newcommand{\GQ}{{\mathcal {GQ}}}
\newcommand{\EE}{{{\mathcal E}(\rho)}}
\newcommand{\HHH}{{\mathbb H}^3}
\newcommand{\CmodTwoPiIZ}{{\mathbf C}/2\pi i {\mathbf Z}}

\def\square{\hfill${\vcenter{\vbox{\hrule height.4pt \hbox{\vrule width.4pt
height7pt \kern7pt \vrule width.4pt} \hrule height.4pt}}}$}

\newenvironment{pf}{\noindent {\it Proof.}\quad}{\square \vskip 6pt}

\title[Length series identities]{A survey of length series identities for surfaces, %
3-manifolds and representation varieties }
\author[Tan, Wong \& Zhang]{Ser Peow Tan, Yan Loi Wong \& Ying Zhang}
\address{Department of Mathematics \\ National University of Singapore \\
2 Science Drive 2 \\ Singapore 117543} \email{mattansp@nus.edu.sg}
\address{Department of Mathematics \\ National University of Singapore \\
2 Science Drive 2 \\ Singapore 117543} \email{matwyl@nus.edu.sg}
\address{Department of Mathematics \\Yangzhou University \\Yangzhou \\ P. R. China  225002
\& IMPA, Estrada Dona Castorina 110, Rio de Janeiro, Brasil
22460-320} \email{yingzhang@yzu.edu.cn; yiing@impa.br}

\thanks{
The authors are partially supported by the National University of
Singapore academic research grant R-146-000-056-112. The third
author is also partially supported by the National Key Basic
Research Fund (China) G1999075104 and a CNPq-TWAS postdoctoral
fellowship.}

%
%

\begin{abstract}
We survey  some of our recent results on length series identities
for hyperbolic (cone) surfaces, possibly with cusps and/or boundary
geodesics; classical Schottky groups; representations/characters of
the one-holed torus group to $\SLTwoC$; and hyperbolic 3 manifolds
obtained by hyperbolic Dehn surgery on punctured torus bundles over
the circle. These can be regarded as generalizations and variations
of McShane's identity for cusped hyperbolic surfaces, which has
found some striking applications in the recent work of Mirzakhani.
We discuss some of the methods and techniques used to obtain these
identities.
\end{abstract}

\maketitle

\section{{\bf Introduction}}\label{s:intro}

In his thesis \cite{mcshane1991thesis}, Greg McShane gave a
remarkable series identity for the lengths of simple closed
geodesics on a complete hyperbolic torus with one cusp. He
generalized this later in \cite{mcshane1998im} to an identity for a
complete hyperbolic surface $M$ with cusps. The identity he obtained
is as follows:

\begin{thm}\label{thm:mcshane general} {\rm (McShane
\cite{mcshane1998im})} In a finite area hyperbolic surface $M$
with cusps and without boundary, let $\Delta_0$ be a distinguished
cusp of $M$. Then
\begin{eqnarray}\label{eqn:McShane}
\sum \frac{1}{1 + \exp \frac{1}{2}(|\alpha|+|\beta|)} =
\frac{1}{2},
\end{eqnarray}
where the sum is taken over all unordered pairs of simple closed
geodesics $\alpha, \beta$ {\rm(}where $\alpha$ or $\beta$ might be
a cusp treated as a simple closed geodesic of length $0${\rm)} on
$M$ such that $\alpha, \beta$ and $\Delta_0$ bound an embedded
pair of pants on $M$, and $|\alpha|$ denotes the length of
$\alpha$.
\end{thm}

In the case of the cusped torus, $\alpha=\beta$ for all pairs in
the sum, and the sum is over all simple closed geodesics $\alpha$
on the torus, which was the original identity obtained in his
thesis.

The proof of the identity was mostly geometric/topological. For
simplicity, consider the case where the surface $M$ has only one
cusp $\Delta_0$. Let ${\mathcal H}$ be the set of all geodesics
emanating from $\Delta_0$. The subset ${\mathcal S}\subset {\mathcal
H}$ of {\em simple} geodesics (no self-intersection) emanating from
$\Delta_0$ turns out to be rather sparse, in fact, by the
Birman-Series Theorem \cite{birman-series1985t}, this set $\mathcal
S$ has zero measure in the set $\mathcal H$. Furthermore, apart from
a countable set of isolated points corresponding to simple geodesics
which also terminate at $\Delta_0$, this set forms a Cantor subset
of $\mathcal H$. Now identifying $\mathcal H$ with a horocycle of
length one about the cusp, it turns out each gap formed by the
complement of the Cantor set has end points corresponding to simple
geodesics which spiral around {\em simple closed geodesics} $\alpha$
and $\beta$ on the surface (with `opposite' spiralling orientation),
where the pair $\alpha$ and $\beta$ bound together with $\Delta_0$
an embedded pair of pants on the surface $M$. Conversely, to every
such pair $\alpha$, $\beta$, there are two gaps with end points
corresponding to simple geodesics which spiral around $\alpha$ and
$\beta$ with opposite orientation, and they both have the same
width. Furthermore, by a simple hyperbolic geometry calculation, the
width of each gap depends only on $|\alpha|$ and $|\beta|$ and is
given by the summand in the left hand side of (\ref{eqn:McShane}).
Theorem \ref{thm:mcshane general} then follows. It should be noted
that the hyperbolic geometry needed to obtain the formula can be
restricted to pairs of pants.

The isolated simple geodesics in ${\mathcal H}$ which start and
end at $\Delta_0$ also have an important geometric interpretation,
each such geodesic $\delta$ defines uniquely a pair of geodesics
$\alpha$ and $\beta$ on $M$ bounding with $\Delta_0$ an embedded
pair of pants in $M$ such that the geodesic $\delta$ is embedded
in the pair of pants. Furthermore,  the two ends of $\delta$ lie
in the corresponding two pairs of gaps.

\vskip 3pt

In brief, the key ingredients in the proof of Theorem
\ref{thm:mcshane general} are:
\begin{itemize}
    \item the study of the set of
simple geodesics on $M$ emanating from the cusp $\Delta_0$,
    \item the Birman-Series theorem, and
    \item some simple geometric identities for
hyperbolic pairs of pants.
\end{itemize}

\vskip 6pt

 There have been several generalizations and variations
of the identity:
\begin{itemize}
    \item  Bowditch gave an independent proof of the
identity for the cusped torus in \cite{bowditch1996blms} and
\cite{bowditch1998plms}, with substantial generalizations to {\em
type-preserving} representations of the punctured torus group to
$\SLTwoC$ satisfying certain conditions, and also a variation of the
identity for complete hyperbolic 3-manifolds which are punctured
torus bundles over the circle in \cite{bowditch1997t}.
    \item  Akiyoshi, Miyachi and Sakuma
gave variations of the identity for quasi-fuchsian punctured torus
groups (in particular, to certain points on the boundary of
quasi-fuchsian space) and hyperbolic punctured surface bundles
over the circle in \cite{akiyoshi-miyachi-sakuma2004cm355} and
\cite{akiyoshi-miyachi-sakuma2004preprint}.
    \item  McShane himself gave variations of the identity arising from a similar analysis of
simple geodesics passing through the Weierstrass points of a cusped
torus and a closed hyperbolic genus two surface in
\cite{mcshane2004blms} and \cite{mcshane2006aasfm}.
    \item  More recently, Mirzakhani proved and
used a version of the identity for bordered hyperbolic surfaces
(surfaces with totally geodesic boundary) to obtain some striking
applications and connections to the Weil-Petersson volume of the
moduli space of bordered Riemann surfaces, the asymptotic behavior
of the number of simple closed geodesics of length less than $L>0$
on a closed hyperbolic surface and the Kontsevich-Witten formula on
the intersection numbers of tautological classes on the moduli space
of curves in \cite{mirzakhani2004preprint}, \cite{mirzakhani2006am}
and \cite{mirzakhani2006jams}. An important observation used in her
papers is that the identity is independent of the hyperbolic
structure on $M$ (with fixed boundary lengths), that is, it holds
for all points in the moduli space.
\end{itemize}

In a different direction, we have also given generalizations and
variations of the identity to
\begin{itemize}
    \item hyperbolic cone surfaces with cusps and/or geodesic
boundary (with all cone angles bounded above by $\pi$), with
applications to generalizations of the Weierstrass identities for
the one-hole/cone torus and closed genus two surface
\cite{tan-wong-zhang2004cone-surfaces};
    \item classical Schottky groups with applications to
    hyperbolic surfaces with geodesic boundary in
\cite{tan-wong-zhang2004schottky}; and
    \item  general (not necessarily
type-preserving) representations of the punctured torus groups to
$\SLTwoC$ with applications to {\em closed} hyperbolic 3-manifolds
obtained by hyperbolic Dehn surgery on hyperbolic punctured torus
bundles over the circle \cite{tan-wong-zhang2004gMm},
\cite{tan-wong-zhang2004necsuf}.
\end{itemize}

In this paper, we will give an exposition of some of the main
results and ideas in \cite{tan-wong-zhang2004cone-surfaces},
\cite{tan-wong-zhang2004schottky}, \cite{tan-wong-zhang2004necsuf}
and \cite{tan-wong-zhang2004gMm}, and also the connection with the
works of McShane, Bowditch, Mirzakhani, and Akiyoshi-Miyachi-Sakuma.

\vskip 3pt

The main point of departure from the works of McShane, Bowditch, and
Akiyoshi-Miyachi-Sakuma in \cite{tan-wong-zhang2004cone-surfaces} is
that we allow $\Delta_0$ (more generally $\Delta_j$) to vary, so
that it is not necessarily a cusp, but may be a cone point or
boundary geodesic. In particular, we consider geometric structures
whose holonomy groups are not necessarily discrete (which represents
also a departure from the point of view of Mirzakhani). From this
point of view, it is natural to consider cone points to have purely
imaginary length. Extending this idea further, more generally, we
may consider representations of the surface group to $\SLTwoC$ so
that the lengths, in particular, the boundary lengths, are not
necessarily real or purely imaginary. From this, we obtain
generalizations of the identity to classical Schottky groups in
\cite{tan-wong-zhang2004schottky}. The main tools are analytic
continuation, a lifting argument, an application of the
Birman-Series argument, and some general comparison results for the
combinatorial length and complex length of an element of the
fundamental group corresponding to a simple closed curve on a marked
surface. In particular, this produces some surprising new identities
for hyperbolic surfaces with boundary, arising from non-standard
markings, for example, we have a nontrivial series identity for the
hyperbolic pair of pants.

Finally, in \cite{tan-wong-zhang2004gMm} and
\cite{tan-wong-zhang2004necsuf}, we pick up on the powerful ideas
and combinatorial techniques of Bowditch to show that a very general
version of the identity can be proved to hold for general $\SLTwoC$
characters of a one-holed torus satisfying some simple conditions.
Similarly, we also show that various relative and restricted
versions of the identity hold. In particular, we are able to give
necessary and sufficient conditions (extended Bowditch Q-conditions)
for the identity to hold for a general $\SLTwoC$ character of the
one-holed torus, \cite{tan-wong-zhang2004necsuf}, and furthermore,
to give relative versions of the identity to characters which are
stabilized by certain cyclic subgroups of the mapping class group
generated either by an Anosov or reducible element, and which
satisfy a relative version of the Bowditch Q-conditions
\cite{tan-wong-zhang2004gMm}. These in turn have applications to
complete and incomplete hyperbolic structures on punctured torus
bundles over the circle, and in particular, give length series
identities for ``almost all'' closed hyperbolic 3-manifolds obtained
by hyperbolic Dehn surgery on a complete hyperbolic torus bundle
over the circle.

A feature of these techniques is that we do not have to use analytic
continuation to obtain the identities, and also, the identities can
be proven for very general representations for which the geometric
interpretation is not necessarily clear. Another interesting feature
of this method is that it gives an independent proof of the
Birman-Series result that the set of simple complete geodesics is
sparse, the point is that in the proof of the series identity, one
is able to prove not just the absolute convergence of the series,
but also to show that a suitably interpreted error term approaches
0.

\vskip 6pt

The rest of this paper is organized as follows. In \S \ref{s:cone},
we discuss the identities for cone surfaces, and also applications
via covering arguments to generalized Weierstrass identities for the
one-holed torus and genus two surface. In \S \ref{s:schottky}, we
discuss the identities for the classical Schottky groups and finally
in \S \ref{s:character}, we discuss the identities for $\SLTwoC$
characters of a one-holed torus.

\vskip 10pt

\noindent {\it Acknowledgements.} The first named author would
like to thank Prof. Michihiko Fujii, the organizer of the
symposium ``Complex Analysis and geometry of hyperbolic spaces"
held at RIMS, Kyoto in Dec 2005 for the invitation to attend and
speak at the symposium. This survey is based on the talks given by
him at the symposium.

\section{{\bf Hyperbolic cone surfaces}}\label{s:cone}
 McShane's original identity (\ref{eqn:McShane}) can be
generalized to hyperbolic cone surfaces, possibly with cusps and/or
totally geodesic boundary, where all cone points have cone angles
less than or equal to $\pi$. For this purpose, it is convenient to
consider the cone points, cusps and boundary geodesics as {\it
geometric} boundary components of $M$ and to define the complex
length of a cone point as $i\theta$, where $\theta$ is the cone
angle, the complex length of a cusp as $0$, and the complex length
of a boundary geodesic as just the usual hyperbolic length. We call
such a surface $M$ a {\em compact hyperbolic cone surface}. We
also define a {\em generalized simple closed geodesic} as %
\begin{enumerate}
\item[(a)] a simple closed geodesic in the geometric interior of $M$; or %
\item[(b)] a geometric boundary component (cone point/cusp/geodesic
boundary) of $M$, with the corresponding complex lengths as
defined earlier; or %
\item[(c)] the double (cover) of a simple geodesic segment joining two
angle $\pi$ cone points on $M$, with length twice the length of the
geodesic segment.
\end{enumerate}
%
The result is then stated as follows.

\begin{thm}\label{thm:complexified} {\rm(Theorem 1.16
\cite{tan-wong-zhang2004cone-surfaces})} Let $M$ be a compact
hyperbolic cone surface with all cone angles in $(0, \pi]$, and
geometric boundary components $\Delta_0, \Delta_1, \cdots, \Delta_N$
with complex lengths $L_0, L_1, \cdots, L_N$ respectively. Then
\begin{eqnarray}\label{eqn:reform of cp and gb cases}
& & \hspace{-20pt}\sum_{\alpha, \,\beta} 2 \tanh^{-1} %
\left( \frac {\sinh\frac{L_0}{2}}{\cosh \frac{L_0}{2}+\exp\frac{|\alpha|+|\beta|}{2}} \right) \nonumber \\ %
& & + \,\,\sum_{j=1}^{N}\sum_{\beta} \tanh^{-1} \left( \frac{\sinh\frac{L_0}{2}\sinh\frac{L_j}{2}} %
{\cosh\frac{|\beta|}{2}+\cosh\frac{L_0}{2}\cosh\frac{L_j}{2}} \right) = \frac{L_0}{2}, %
\end{eqnarray}
if $\Delta_0$ is a cone point or a boundary geodesic; and
\begin{eqnarray}\label{eqn:reform of cusp cases}
\sum_{\alpha,\,\beta}\frac{1}{1+\exp\frac{|\alpha|+|\beta|}{2}}+
\sum_{j=1}^{N}\sum_{\beta}\frac{1}{2}\frac{\sinh\frac{L_j}{2}}
{\cosh\frac{|\beta|}{2}+\cosh\frac{L_j}{2}}=\frac{1}{2},
\end{eqnarray}
if $\Delta_0$ is a cusp; where in either case the first sum is taken
over all unordered pairs of generalized simple closed geodesics
$\alpha, \beta$ on $M$ which bound with $\Delta_0$ an embedded pair
of pants on $M$ {\rm(}note that one of $\alpha, \beta$ might be a
geometric boundary component\,{\rm)} and the sub-sum in the second
sum is taken over all generalized simple closed geodesics $\beta$
which bounds with $\Delta_j$ and $\Delta_0$ an embedded pair of
pants on $M$.
Furthermore, each series in {\rm(\ref{eqn:reform of cp and gb
cases})} and {\rm(\ref{eqn:reform of cusp cases})} converges
absolutely.
\end{thm}


The summands in the first sum correspond to main gaps and those in
the second series correspond to side gaps. Note that McShane's
identity (\ref{eqn:McShane}) is a special case of (\ref{eqn:reform
of cusp cases}) where $\Delta_0$ is a cusp and all the summands in
the second series are zero since none of $\Delta_j$, $j=1, \ldots N$
are cone points or boundary geodesics. Also, the identity
(\ref{eqn:reform of cusp cases}) can be derived from the first order
infinitesimal terms of the identity (\ref{eqn:reform of cp and gb
cases}). \vskip 5pt

For the purpose of generalizations to classical Schottky groups
later, we define the functions $G(x,y,z)$ and $S(x,y,z)$
corresponding to the ``main gaps'' and the ``side gaps'' as follows,
where the $\log$ function takes its principal branch, i.e., with
imaginary part in $(-\pi, \pi]$, and the function $\tanh^{-1}$ is
defined by %
$$\tanh^{-1}(x)=\frac{1}{2}\log\frac{1+x}{1-x} \quad \quad %
\textrm{for} \,\, x \in \mathbf C \backslash \{\pm 1\},$$ %
and hence has imaginary part in $(-\pi/2, \pi/2)$.

\begin{defn}\label{def:GandS} For $x,y,z \in \mathbf C$, we define
\begin{eqnarray}
G(x,y,z)\!\!&:=&\!\!2\tanh^{-1}\left(\frac{\sinh(x)}{\cosh(x)+\exp(y+z)}\right), \\ %
S(x,y,z)\!\!&:=&\!\!\tanh^{-1}\left(\frac{\sinh(x)\sinh(y)}{\cosh(z)+\cosh(x)\cosh(y)}\right).
\end{eqnarray}
\end{defn}

\noindent It can be shown that $G(x,y,z)$ and $S(x,y,z)$ can also be
expressed as
\begin{eqnarray}
G(x,y,z)\!\!&=&\!\!\log\frac{\exp(x)+\exp(y+z)}{\exp(-x)+\exp(y+z)}, \label{eqn:G in log} \\ %
S(x,y,z)\!\!&=&\!\!\frac{1}{2}\log\frac{\cosh(z)+\cosh(x+y)}{\cosh(z)+\cosh(x-y)}, \label{eqn:S in log} %
\end{eqnarray}

\noindent as used by Mirzakhani in \cite{mirzakhani2004preprint}.

\vskip 6pt

 The basic idea of the proof of Theorem \ref{thm:complexified} is similar to that in
\cite{mcshane1998im}, we pick a distinguished boundary component
$\Delta_0$ (which may be a cone point, cusp, or boundary geodesic),
and consider the set of all geodesics ${\mathcal H}$ emanating
normally from $\Delta_0$ (one only needs to worry about ``normally''
when $\Delta_0$ is a boundary geodesic). Topologically, this set is
a circle; geometrically, we can put a natural measure on this circle
as follows. In the case $\Delta_0$ is a cusp, we identify ${\mathcal
H}$ as before with the horocycle of length one around $\Delta_0$; in
the case $\Delta_0$ is a boundary geodesic of length $L_0$, we
identify ${\mathcal H}$ with $\Delta_0$ itself with length $L_0$;
and in the case where $\Delta_0$ is a cone point with cone angle
$\theta_0$, we identify ${\mathcal H}$ with a circle about
$\Delta_0$ with the natural radian measure $\theta_0$. We then
consider the subset ${\mathcal S} \subset {\mathcal H}$ consisting
of the {\em simple, complete} geodesics, by which, we mean the
geodesics emanating normally from $\Delta_0$ which are simple and do
not terminate at a cone point, cusp, or boundary geodesic, hence are
complete in the forward direction. By a slight variation of the
Birman-Series Theorem, this set again has zero measure in ${\mathcal
H}$, and, as in the previous case, is essentially a Cantor set
(there may be a countable collection of isolated points). As before,
the complement ${\mathcal H}\setminus {\mathcal S}$ consists of gaps
bounded by end points which correspond to geodesics spiralling
around simple closed curves. However, in this case, besides the main
gaps which we had before, side gaps can occur, if some of the other
boundary components $\Delta_k$ are boundary geodesics. In this case,
a side gap is bounded by two points corresponding to simple
geodesics that spiral around the same boundary component, but in
opposite directions.

\vskip 3pt

 It turns out that the analysis is actually easier if
we look at the set of geodesics in ${\mathcal H}$ which are not in
${\mathcal S}$, that is, are not simple and complete. In this
case, the geodesic either intersects a boundary component, or has
a self intersection, in the latter case, we consider the initial
part of the geodesic up to the first point of self intersection.
In either case, by considering a tubular neighborhood of the union
of this geodesic (segment) $\delta$ with $\Delta_0$, one obtains
two curves $\alpha$ and $\beta$, unique up to homotopy which bound
together with $\Delta_0$ an embedded pair of pants $P$ in $M$
which contains $\delta$. Now the condition that all cone angles
are less than or equal to $\pi$ ensures that $\alpha$ and $\beta$
are realizable as generalized geodesics, so we have an embedded
pair of pants with $\Delta_0$, $\alpha$ and $\beta$ as the
boundary components.  Now the geodesic $\delta$ lies in either a
main gap or a side gap, see Figures \ref{fig01} and \ref{fig02},
where the geodesics $\gamma_{\alpha}$ and $\gamma_{\beta}$ in
Figure \ref{fig02} are geodesics which spiral around $\alpha$ and
$\beta$ with opposite orientations and bound a main gap. The
computation of the width of the gaps proceeds as before but is
somewhat more complicated because of the various cases that can
occur depending on whether or not one of $\alpha$, $\beta$ is a
boundary geodesic around a cone point of $M$. Note that if
$\Delta_j$ is a cone point, we do not expect to have a side gap
from this point of view, however, if we wish to interpret the gaps
as analytic functions of the boundary lengths, then there should
be a purely imaginary side gap if for example $\Delta_0$ is a
boundary geodesic and $\Delta_j$ is a cone point. Similarly, in
this case, the main gap is no longer real or purely imaginary. In
fact, the formula given in Theorem \ref{thm:complexified} takes
this analytic point of view and is a complexified, unified version
of all these different cases. There is a geometric interpretation
of these (complexified) gaps by considering the picture in $\HHH$,
see \cite{tan-wong-zhang2004cone-surfaces} for details.

In the case where there are no cone points, $G(x,y,z)$ and
$S(x,y,z)$ are positive real for all summands, and the absolute
convergence is trivial, but if there are some cone points,   the
summands in the formula are not necessarily real and positive, and
the absolute convergence of the various series is no longer
obvious and requires justification, hence the last statement given
in the theorem.  The absolute convergence is proven  by using a
modification of the Birman-Series argument.


 It is important to note that for
the above analysis to work, all essential simple closed curves
should be realizable as generalized simple geodesics, that is,
either geodesics or the double cover of a geodesic segment between
two angle $\pi$ cone points. However, for this to be true, we
require all cone angles to be less than or equal to $ \pi$, and our
proof is by a convexity argument and a suitable application of the
Arzela-Ascoli Theorem. It is not clear how this condition can be
relaxed, hence, a McShane type identity for general closed
hyperbolic surfaces without boundary remains elusive.


 The above is closely related to the formula obtained by Mirzakhani
for hyperbolic surfaces with geodesic boundary in
\cite{mirzakhani2004preprint}. In particular, her analysis works for
the cone hyperbolic surfaces we consider and the same (recursive)
formula for the Weil-Petersson volumes of the moduli space of
bordered Riemann surfaces holds for cone Riemann surfaces (possibly
with geodesic boundary), where the lengths of cone boundary
components are given by $i\theta$, where $\theta$ is the cone angle.
It also seems that the same analysis she uses to study the
asymptotics of the lengths of simple closed geodesics on closed
hyperbolic surfaces in \cite{mirzakhani2006am} should carry over to
the situation we study, namely, there should be a constant $C_M$
depending only on the hyperbolic structure on the cone hyperbolic
surface $M$ (with all cone angles bounded above by $\pi$) such that
the number of simple closed geodesics on $M$ of length less than $L$
is asymptotic to $C_M\cdot L^{6g-6+2N}$ where $N$ is the number of
geometric boundary components, and $6g-6+2N>0$. As for the relation
to the Kontsevich-Witten formula, and the recursion formula for the
volumes of the moduli space in \cite{mirzakhani2006jams}, there is
also some recent work of Do and Norbury \cite{do-norbury2006arxiv}
generalizing Mirzakhani's work to cone surfaces.

\vskip 3pt

We summarize the various points raised above:
 \begin{itemize}
    \item For a cone hyperbolic surface $M$ possibly with cusps and/or
    geodesic boundary, if all cone angles are less than or equal to $\pi$,
    then all essential simple closed curves on $M$ are realizable by
    (generalized) simple closed geodesics.
    \item The Birman-Series theorem generalizes to these cone surfaces.
    A modification of the argument used in the proof can
    also be used to prove the absolute convergence of the various
    series in the identity.
    \item The gaps formed by taking the complement of the simple complete geodesics
    emanating normally from a fixed boundary component can be
    calculated. Apart from the main gaps which occur in the cusped
    case, side gaps may also occur if there are other boundary
    components which are cone points or boundary geodesics.
    \item It is easier to study the set of geodesics which are not
    simple and complete, these either have self intersection or
    intersect the boundary of $M$, and give rise to pairs of pants
    embedded in the surface.
    \item The analysis of geodesics in $M$ emanating from a boundary
    component can be restricted to just the analysis of geodesics
    in a pair of pants.
 \end{itemize}

\vskip 3pt

Theorem \ref{thm:complexified} together with the fact that a
one-holed hyperbolic torus/(closed hyperbolic surface of genus two)
admits a canonical elliptic/(hyperelliptic) involution and some
general covering arguments can be used to deduce further identities
for the one-holed hyperbolic torus/(genus two surface). These can be
regarded as generalizations of the Weierstrass identities given by
McShane in \cite{mcshane2004blms} and \cite{mcshane2006aasfm}. Here
when we say a one-holed torus, we mean that the boundary may be a
geodesic, cusp or cone point. We have:

\vskip 4pt
\begin{cor}\label{cor:mcshane conical holed weierstrass}{\rm(Corollary
1.10 \cite{tan-wong-zhang2004cone-surfaces})} Let $T$ be either a
hyperbolic one-cone torus where the single cone point has cone angle
$\theta \in [0, 2\pi)$ or a hyperbolic one-holed torus where the
single boundary geodesic has length $l\ge 0$. Then we have
respectively
\begin{eqnarray}
\sum_{\gamma \in \mathcal A} \tan^{-1} \bigg ( \frac{\cos
\frac{\theta}{4}}{\sinh \frac{|\gamma|}{2}} \bigg )=
\frac{\pi}{2},
\end{eqnarray}
\begin{eqnarray}
\sum_{\gamma \in \mathcal A} \tan^{-1} \bigg ( \frac{\cosh
\frac{l}{4}}{\sinh \frac{|\gamma|}{2}} \bigg )= \frac{\pi}{2},
\end{eqnarray}
where the sum in either case is taken over all the simple closed
geodesics $\gamma$ in a given Weierstrass class $\mathcal A$.
\end{cor}

Note that a cusp can be regarded either as a cone point of cone
angle $0$ or a geodesic of length 0 in either of the cases above.

\vskip 4pt

\begin{thm}\label{thm:mcshane genus two global}{\rm(Theorem 1.13,
\cite{tan-wong-zhang2004cone-surfaces})} %
Let $M$ be a genus two closed hyperbolic surface. Then
\begin{eqnarray}
\sum \tan^{-1} \exp \left ( -\frac{|\alpha|}{4}- \frac{|\beta|}{2} \,\right ) %
= \frac{3\pi}{2}, %
\end{eqnarray}
where the sum is taken over all ordered pairs $(\alpha, \beta)$ of
disjoint simple closed geodesics on $M$ such that $\alpha$ is
separating and $\beta$ is non-separating.
\end{thm}

\vskip 4pt

In fact, Corollary \ref{cor:mcshane conical holed weierstrass} can
be extended to much more general representations of $\pi_1(T)$ to
$\SLTwoC$ (see \cite{tan-wong-zhang2004schottky}) and Theorem
\ref{thm:mcshane genus two global} can be extended to
quasi-fuchsian representations of $\pi_1(M)$ to $\PSLTwoC$
(\cite{tan-wong-zhang2004cone-surfaces} Addendum 1.15).

\vskip 6pt

\noindent {\sc Sketch of proof of Corollary \ref{cor:mcshane conical
holed weierstrass} and Theorem \ref{thm:mcshane genus two global}.}
Let $\iota$ be the elliptic involution on $T$. Then $T/\iota$ is a
sphere with four boundary components, three of which are cone points
of angle $\pi$ and the fourth a boundary component of length $l/2$
or a cone point of cone angle $\theta/2$ depending on whether $T$
has a boundary geodesic of length $l$ or a cone point of angle
$\theta$, respectively. Apply Theorem \ref{thm:complexified} to
$T/\iota$ with one of the cone points of angle $\pi$ as $\Delta_0$.
Then the sum is over all generalized simple closed geodesics on
$T/\iota$ which are double covers of geodesic segments joining the
other two cone points of angle $\pi$, these lift to geodesics on $T$
which are in the Weierstrass class consisting of all geodesics which
miss the lift of $\Delta_0$ on $T$, giving Corollary
\ref{cor:mcshane conical holed weierstrass}. For Theorem
\ref{thm:mcshane genus two global}, again consider the hyperelliptic
involution $\iota$ on $M$. Then $M/\iota$ is a sphere with six cone
points, all of cone angle $\pi$. Apply Theorem
\ref{thm:complexified} to each of the six cone points. For each
identity, the sum is now over all pairs of disjoint $\alpha'$ and
$\beta'$ on $M/\iota$ such that $\alpha'$ is a geodesic on $M/\iota$
which separates it to two pieces each containing three cone points,
and $\beta'$ is a double cover of a geodesic segment on the piece
separated by $\alpha'$ containing $\Delta_0$ which connects the
other two cone points. Now take the sum over all the six identities
and lift the result to $M$. Note that $\alpha'$ lifts to a
separating geodesic on $M$ and $\beta'$ lifts to a disjoint
non-separating geodesic on $M$ and furthermore, all separating
geodesics on $M$ project to separating geodesics on $M/\iota$ which
separate $M/\iota$ to two components each containing exactly three
cone points while non-separating geodesics on $M$ project to
geodesic arcs on $M/\iota$ connecting exactly two of the cone
points.


\begin{figure}
\begin{center}
\mbox{\beginpicture \setcoordinatesystem units <0.24in,0.24in> %
\setplotarea x from 0 to 9, y from 0 to 6 %
\ellipticalarc axes ratio 2:1 360 degrees from 3.6 5.0 center at 5.0 5.0 %
\setquadratic
\setdashes<2pt> \plot 9.560660172  2.560660172  9.399493762
2.618001774  9.150278635  2.565905573  8.837409705  2.409471107
8.491512762  2.164011274  8.146446610  1.853553390  7.835988726
1.508487238  7.590528893  1.162590295  7.434094427  0.8497213651
7.381998226  0.6005062381  7.439339828  0.439339828 / %
\plot 2.560660172 0.439339828  2.618001774  0.6005062381 2.565905573
0.8497213651  2.409471107  1.162590295  2.164011274 1.508487238
1.853553390  1.853553390  1.508487238  2.164011274 1.162590295
2.409471107  0.8497213651  2.565905573  0.6005062381
2.618001774 0.439339828  2.560660172 / %
\setsolid \plot 9.560660172  2.560660172 9.618001774  2.399493762
9.565905573  2.150278635  9.409471107 1.837409705  9.164011274
1.491512762  8.853553390  1.146446610 8.508487238  0.8359887260
8.162590295  0.5905288924  7.849721365
0.4340944273  7.600506238  0.3819982259  7.439339828  0.439339828 / %
\plot 2.560660172  0.439339828  2.399493762  0.3819982259
2.150278635 0.4340944273  1.837409705  0.5905288924  1.491512762
0.8359887260 1.146446610  1.146446610  0.8359887260  1.491512762
0.5905288924 1.837409705  0.4340944273  2.150278635  0.3819982259
2.399493762 0.439339828  2.560660172 / %
\plot 9.560660172  2.560660172 7.5 4.0 6.4 5.0 / %
\plot 0.439339828 2.560660172 2.5 4.0 3.6 5.0 / %
\plot 2.560660172  0.439339828 5.0 1.0 7.439339828  0.439339828 / %
\plot 5.224157199 4.3 5.1 1.7 5.0 1.0  / %
\setdashes<2pt> \plot 5.0 1.0 4.8 2.0 4.775842801 5.75 / %
\setsolid \plot 4.628894061 5.66 4.62 5.5 4.76 5.5 / %
\plot 5.39 4.3 5.39 4.15 5.21 4.14 / \plot 4.25 4.4 4.1 1.7 3.75 0.85 / %
\setdashes<2pt> \plot 3.75 0.85 2.4 2.1 2 3.6 / %
\setsolid \plot 2 3.6 3 2.6  3.5 0.78 / %
\setdashes<2pt> \plot 3.5 0.78 2.3 1.5 1.2 3.0 / %
\setsolid \plot 1.2 3.0 1.6 3.1 3.2 2.6 / \plot 4.4 4.35 4.4 4.2 4.25 4.23 / %
\put {\mbox{\LARGE $\cdot$}} [cc] <0mm,0mm> at 2.93 2.66 %
\arrow <6pt> [0.16,0.6] from 5.19 3.7 to 5.188 3.6 %
\arrow <6pt> [0.16,0.6] from  4.23 3.7 to 4.22 3.6 %
\put {\mbox{\small $\triangle_0$}} [cb] <0mm,2mm> at 5.0 5.5 %
\put {\mbox{\small $\alpha$}} [rt] <-1mm,-1.6mm> at 1.4 1.2 %
\put {\mbox{\small $\beta$}} [lt] <1mm,-1mm> at 9.0 1.2 %
\put {\mbox{\small $\gamma$}} [lc] <1mm,0mm> at 5.1 3.2 %
\put {\mbox{\small $\delta$}} [rc] <-1mm,0mm> at 4.3 3.2 %
\endpicture}
\hspace{0.4in} \mbox{\beginpicture
\setcoordinatesystem units <0.29in,0.29in> %
\setplotarea x from 3.5 to 9, y from 0.5 to 6 %
\plot 9.560660172 2.560660172 7 5 4 2  / %
\setquadratic \plot 9.560660172 2.560660172 9.639852575  2.377642961
9.607468267 2.108715941 9.466677248 1.780203564  9.231261125
1.424262911 8.924264069 1.075735931 8.575737089  0.7687388748
8.219796436 0.5333227521 7.891284059 0.3925317336  7.622357039
0.3601474247 7.439339828 0.439339828 / %
\setdashes<2pt> \plot
9.560660172 2.560660172 9.377642961 2.639852575  9.108715941
2.607468267 8.780203564 2.466677248 8.424262911  2.231261125
8.075735931 1.924264069 7.768738875 1.575737089  7.533322752
1.219796435 7.392531733 0.8912840588 7.360147425  0.6223570393
7.439339828 0.439339828 / %
\put {\mbox{\LARGE $\cdot$}} [cc] <0mm,0mm> at 7 5 %
\put {\mbox{\scriptsize $\triangle_{0}$}} [cb] <0mm,1mm> at 7 5 %
\put {\mbox{\LARGE $\cdot$}} [cc] <0mm,0mm> at 4 2 %
\put {\mbox{\scriptsize $\alpha$}} [rt] <-1mm,0mm> at 4 2 %
\put {\mbox{\scriptsize $\beta$}} [rt] <40mm,0mm> at 4 2 %
\setsolid \plot 4.0 2.0 6.3 1.2 7.439339828  0.439339828 / %
\setdashes<2pt> \plot 5.5 1.55 5.1 2.0 5.0 3.0 / %
\setsolid \plot 5 3 5.3 2.3 5.2 1.62 / %
\setdashes<2pt> \plot 5.2 1.62 4.8 2  4.7 2.7 / %
\setsolid \plot 4.7 2.7 5 2.7 5.5 2.3 / %
\put {\mbox{\LARGE $\cdot$}} [cc] <0mm,0mm> at 5.25 2.5 %
\setsolid \plot 7 5 6.1 2.3 5.5 1.55 / \plot 6.6 1.25 7.1 2.5 7 5 / %
\setdashes<2pt> \plot 7 5 6.3 1.8 6.6 1.25 / %
\setsolid \arrow <6pt> [.16,.6] from  6.19 2.5 to 6.15 2.4 %
\arrow <6pt> [0.16,0.6] from  7.02 2 to 7 1.9 %
\put {\mbox{\scriptsize $\gamma$}} [lt] <0.6mm,0mm> at 7 2 %
\put {\mbox{\scriptsize $\delta$}} [rc] <-0.6mm,0mm> at 6.2 2.5 %
\endpicture}
\end{center}
\caption{}\label{fig01}
\end{figure}


\vskip 4pt

\begin{figure}
\begin{center}
\hskip -20pt \mbox{\beginpicture \setcoordinatesystem units <0.225in,0.225in> %
\setplotarea x from 1 to 9, y from 0 to 6 %
\ellipticalarc axes ratio 2:1 360 degrees from 3.6 5 center at 5 5 %
\setquadratic \setdashes<2pt> \plot 9.560660172 2.560660172
9.344866759 2.672628777  9.046371901 2.669812307 8.694394354
2.552486458 8.323388134  2.332135902 7.969669914 2.030330086
7.667864098 1.676611866  7.447513542 1.305605646 7.330187693
0.9536280995 7.327371223  0.6551332412 7.439339828 0.439339828 / %
\setsolid \plot 9.560660172  2.560660172  9.672628777 2.344866759
9.669812307 2.046371901  9.552486458  1.694394354 9.332135902
1.323388134 9.030330086  0.9696699142  8.676611866 0.6678640980
8.305605646 0.4475135417  7.953628099  0.3301876929 7.655133241
0.3273712228 7.439339828  0.439339828 / %
\plot 2.560660172 0.439339828 2.669812307  0.9536280995  2.332135902
1.676611866 1.676611866 2.332135902  0.9536280995  2.669812307
0.439339828 2.560660172 0.3301876929  2.046371901  0.6678640980
1.323388134 1.323388134 0.6678640980  2.046371901  0.3301876929
2.560660172 0.439339828 / %
\plot 9.560660172  2.560660172 7.5 4.0 6.4 5.0 / %
\plot 0.439339828  2.560660172 2.5 4.0 3.6 5.0 / %
\plot 2.560660172  0.439339828 5.0 1.0 7.439339828  0.439339828 / %
\plot 5.224157199 4.3 5.1 1.7 5.0 1.0  / %
\setdashes<2pt> \plot 5.0 1.0 4.8 2.0  4.775842801 5.7 / %
\setlinear\setsolid \plot 4.6 5.66 4.62 5.5 4.76 5.5 / %
\plot 5.39 4.3 5.39 4.15 5.21 4.14 / %
\setquadratic \plot 5.75 4.4 5.9 1.7 6.0 0.9 / %
\setlinear \plot 5.9 4.45 5.9 4.3 5.8 4.27 / %
\setquadratic \setdashes<2pt> \plot 6.0 0.9 7.0 3.0 8.0 3.6 / %
\setsolid \plot 8.0 3.6 8.5 2.8 7.8 1.4 /
\plot 4.25 4.4 4.1 1.7 3.75 0.85 / %
\setdashes<2pt> \plot 3.75 0.85 2.4 2.1 2.0 3.6 / %
\setsolid \plot 2.0 3.6 3.0 2.6  3.5 0.78 / %
\setdashes<2pt> \plot 3.5 0.78 3.0 1.0 2.5 1.5 / %
\setlinear \setsolid \plot 4.4 4.35 4.4 4.2 4.25 4.23 / %
\arrow <6pt> [.16,.6] from 5.19 3.7 to 5.188 3.6 %
\arrow <6pt> [0.16,0.6] from  4.23 3.7 to 4.22 3.6 %
\arrow <6pt> [0.16,0.6] from  5.79 3.7 to 5.8 3.6 %
\put {\mbox{\small $\triangle_0$}} [cb] <0mm,2mm> at 5.0 5.5 %
\put {\mbox{\small $\delta_\alpha$}} [rb] <-1mm,1mm> at 3.0 4.2 %
\put {\mbox{\small $\delta_\beta$}} [lb] <1mm,1mm> at 7.0 4.2 %
\put {\mbox{\small $\alpha$}} [rt] <-1mm,-1.6mm> at 1.0 1.2 %
\put {\mbox{\small $\beta$}} [lt] <1mm,-1mm> at 9.0 1.2 %
\put {\mbox{\small $\gamma$}} [lc] <1mm,0mm> at 5.1 3.2 %
\put {\mbox{\small $\gamma_{{}_{\alpha}}$}} [rc] <-1mm,0mm> at 4.3 3.2 %
\put {\mbox{\small $\gamma_{{}_{\beta}}$}} [lc] <1.5mm,0mm> at 5.7 3.2 %
\endpicture}
\hspace{0.4in}
\mbox{\beginpicture \setcoordinatesystem units <0.27in,0.27in> %
\setplotarea x from 3.5 to 9.0, y from 0.2 to 6.0 %
\plot 9.560660172 2.560660172 7.0 5.0 4.0 2.0  / %
\setquadratic \plot 9.560660172 2.560660172 9.639852575 2.377642961
9.607468267 2.108715941 9.466677248 1.780203564 9.231261125
1.424262911 8.924264069 1.075735931 8.575737089 0.7687388748
8.219796436 0.5333227521 7.891284059 0.3925317336 7.622357039
0.3601474247 7.439339828 0.439339828 / %
\setdashes<2pt> \plot
9.560660172 2.560660172 9.377642961 2.639852575  9.108715941
2.607468267 8.780203564 2.466677248 8.424262911  2.231261125
8.075735931 1.924264069 7.768738875 1.575737089  7.533322752
1.219796435 7.392531733
0.8912840588 7.360147425  0.6223570393 7.439339828 0.439339828 / %
\put {\mbox{\Huge $\cdot$}} [cc] <0mm,0mm> at 7 5 %
\put {\mbox{\small $\triangle_{0}$}} [cb] <0mm,1mm> at 7 5 %
\put {\mbox{\LARGE $\cdot$}} [cc] <0mm,0mm> at 4 2 %
\put {\mbox{\small $\alpha$}} [rt] <-1mm,0mm> at 4 2 %
\put {\mbox{\small $\beta$}} [rt] <0mm,0mm> at 9.2 1 %
\setsolid \plot 4 2 6.3 1.2 7.439339828  0.439339828 / %
\setdashes<2pt> \plot 7.0 5.0 5.9 2.9 5.5 1.55 / %
\setsolid \plot 7 5 6.1 2.3 5.5 1.55 / \plot 7 5 6.5 1.8 6.6 1.25 / %
\setdashes<2pt> \plot 6.6 1.25 8.4 2.9 8.9 3.1 / %
\setsolid \plot 8.9 3.1 8.7 2.5 8.2 1.6 / %
\arrow <6pt> [0.16,0.6] from  6.19 2.5 to 6.15 2.4 %
\arrow <6pt> [0.16,0.6] from 6.51 2.0 to 6.49 1.9 %
\put {\mbox{\small $\gamma_{\beta}$}} [lt] <2.0mm,0mm> at 6.3 2.0 %
\put {\mbox{\small $\gamma$}} [rc] <-0.6mm,0mm> at 6.2 2.5 %
\put {\mbox{\small $\delta_\beta$}} [lb] <0mm,0mm> at 8.0 4.0 %
\put {\mbox{\small $\delta_{\alpha}\hspace{-0.03in}=\hspace{-0.02in}\gamma_{\alpha}$}} %
[rc] <-1.5mm,0mm> at 6.2 4.2 %
\endpicture}
\end{center}
\caption{}\label{fig02}
\end{figure}

\section{{\bf Classical Schottky groups}}\label{s:schottky}

We first note that if $M$ is a hyperbolic surface with geodesic
boundary components, then the holonomy group is in fact a fuchsian
Schottky group. We next observe that in Theorem
\ref{thm:complexified}, the summands in the series are all
analytic functions of the lengths (if we take the analytic
continuation of the $\tanh^{-1}$ function). These are (real)
analytic in the parameters of the Teichm\"uller space, which in
turn is locally homeomorphic to the representation variety (modulo
conjugation) of representations from $\pi_1(M)$ to $\PSLTwoR$. It
is natural to see if we can apply analytic continuation to obtain
generalizations of the result to representations of $\pi_1(M)$ to
$\PSLTwoC$ or $\SLTwoC$. The absolute convergence of the series in
question and the connectedness of the deformation space are the
two key issues. There are other important technicalities. It turns
out we can do this and obtain series identities for classical
Schottky groups which generalize McShane's identity. We summarize
below some of the relevant points that crop up:
\begin{itemize}
    \item Absolute convergence of the series in
    (\ref{eqn:reform of cp and gb cases}) for classical Schottky groups;
    \item Connectedness of the deformation space;
    \item Lifting of the representations from $\PSLTwoC$ to
    $\SLTwoC$;
    \item Determination of an explicit half-length for transformations
    in $\SLTwoC$;
    \item Choice of a fuchsian marking that will determine how
    the summands in (\ref{eqn:reform of cp and gb cases}) are obtained.

\end{itemize}

To start with, we define classical Schottky space. Fix $n \ge 2$.
This is the space of (marked) faithful representations from the free
group $F_n$ on $n$ generators to $\PSLTwoC$, up to conjugation, such
that the image is a classical Schottky group. We keep track of the
marking, as this makes the statement of the results clearer and more
precise later.

\begin{defn}\label{def:markedschottkygroup}{\rm
A {\it {\rm(}marked\,{\rm)} classical Schottky group} (of rank $n$)
is a discrete, faithful representation $\rho:F_n \to \PSLTwoC$ such
that there is a region $D \subset \Cinfty $, where $D$ is bounded by
$2n$ disjoint geometric circles  $C_1, C'_1, \cdots, C_n, C'_n$ in
$\Cinfty$, so that, for $i=1, \ldots, n$, $\rho(a_i)(C_i)= C'_i$,
and $\rho(a_i)(D) \cap D = \emptyset$. It is {\it fuchsian} if the
representation can be conjugated to a representation into
$\PSLTwoR$.} Two representations are equivalent if they are
conjugate by an element of $\PSLTwoC$.
\end{defn}

The space of equivalent classes of marked classical Schottky groups
is the marked classical Schottky space, denoted by ${\mathcal
S}^{\rm mc}_{\rm alg}$. To simplify notation, we use $\rho$ instead
of $[\rho]$ to denote the elements of ${\mathcal S}^{\rm mc}_{\rm
alg}$. Note that every element of a classical Schottky group is
loxodromic. One may associate a complex length $l(A)$ to each
loxodromic element $A \in \PSLTwoC$, where if we consider $A$ as an
orientation preserving isometry of $\HHH$, the real part of $l(A)$
is the (positive) translation distance of $A$ along its axis, and
the imaginary part is the rotation about the axis, where the
orientation is naturally induced by the translation direction of
$A$. The complex length $l(A)$ is related to the trace by the
formula
\begin{equation}\label{eqn:complexlength}
    l(A)=2\cosh^{-1}\!\big(\!\textstyle-\frac12{\rm tr}\,(A)\big),
\end{equation}
and is chosen to have positive real part (note that we could have
done away with the minus sign inside the $\cosh^{-1}$ function since
the trace is only defined up to $\pm$ sign, we add it here for
consistency with the definition for the half length to be given
later). Then $l(A)$ is defined up to multiples of $2 \pi i$, and
depend only on $\pm {\rm tr}\,(A)$ or ${\rm tr}^2\,(A)$. More
explicitly, we have $l(A)=\cosh^{-1}\!(\frac12{\rm tr}^2(A)-1)$.

We may give a natural parametrization of ${\mathcal S}^{\rm
mc}_{\rm alg}$ by the ideal fixed points, and the square of the
traces or the complex lengths of $\rho(a_i)$, $i=1, \ldots, n$ as
follows; here we use ${\rm Fix}^{\pm}\rho(a_i)$ to denote the
attracting and repelling fixed points of $\rho(a_i)$.

We first normalize $\rho$  by conjugation so that
\vskip 3pt %
\centerline{${\rm Fix}^{-}\rho(a_1)=0, \quad {\rm Fix}^{+}\rho(a_1)=\infty$ \,\, and \,\, ${\rm Fix}^{-}\rho(a_2)=1$.} %
\vskip 3pt %
\nnn Then it is not difficult to see that we can parameterize $\rho$ by %
\begin{eqnarray*}
&&(\,{\rm Fix}^{+}\rho(a_2), {\rm Fix}^{-}\rho(a_3), {\rm Fix}^{+}\rho(a_3), %
\cdots, {\rm Fix}^{+}\rho(a_n); \, {\rm tr}^2\, \rho(a_1),\cdots, {\rm tr}^2\,\rho( a_n)) \\ %
&&\in \Cinfty^{2n-3} \times {\mathbf C}^{n}, %
\end{eqnarray*}
or, alternatively, by
\begin{eqnarray*}
&&(\,{\rm Fix}^{+}\rho(a_2), {\rm Fix}^{-}\rho(a_3), {\rm Fix}^{+}\rho(a_3), %
\cdots, {\rm Fix}^{+}\rho(a_n); \,l(\rho(a_1)),\cdots, l(\rho( a_n))) \\ %
&&\in \Cinfty^{2n-3} \times (\CmodTwoPiIZ)^{n}. %
\end{eqnarray*}

With this normalized parametrization we have

\begin{lem}\label{lem:connected} {\rm (Maskit \cite{maskit2002cm311})}
The marked classical Schottky space ${\mathcal S}^{\rm mc}_{\rm
alg}$ is a path connected open subset of $\Cinfty^{2n-3} \times
(\CmodTwoPiIZ)^{n}$.
\end{lem}

\begin{defn}\label{def:fuchsianmarking}{\rm
A fuchsian marking in ${\mathcal S}^{\rm mc}_{\rm alg}$ is a {\em
fuchsian} representation $\rho_0 \in {\mathcal S}^{\rm mc}_{\rm
alg}$.}
\end{defn}

For a fuchsian marking $\rho_0$, $\HH/\rho_0(F_n)$ is a complete
hyperbolic surface. Its convex core, $M_0$, is a hyperbolic
surface with geodesic boundary, which we call the hyperbolic
surface corresponding to the fuchsian marking. Let $\Delta_0,
\Delta_1, \ldots, \Delta_m$ be the boundary components of $M_0$.
The image $\rho_0(F_n)$, and hence $F_n$ (since $\rho_0$ is
faithful), can be identified with $\pi_1(M_0)$, and if we define
an equivalence relation $\sim$ on $F_n$ by $g \sim h$ if $g$ is
conjugate to $h$ or $h^{-1}$, then there is a bijection
$${\mathfrak f}:F_n/\sim ~ \to {\mathcal C}$$ from $F_n/\sim$ to the set
${\mathcal C}$ of free homotopy classes of closed curves on $M_0$.
Note that there is a unique geodesic representative on $M_0$ for
each nontrivial element of ${\mathcal C}$.

\begin{defn}\label{def:sets}
For a fixed fuchsian marking $\rho_0$, let $M_0$ be the
corresponding hyperbolic surface. Let $\Delta_0, \Delta_1, \ldots,
\Delta_n$ be the boundary components of $M_0$, and let $[d_i]\in
F_n/\sim$, $i=0,\ldots,m$ be the equivalence class corresponding
to the boundary component $\Delta_i$, that is ${\mathfrak
f}[d_i]=\Delta_i$.
\begin{itemize}
    \item [(a)] We define ${\mathcal P}$ to be the set of
all unordered pairs $\{[g],[h]\}$ of elements in $F_n/\sim$ such
that ${\mathfrak f}[g]$ and ${\mathfrak f}[h]$ are free homotopy
classes of simple closed curves which bound together with $\Delta_0$
an embedded pair of pants in $M_0$ (note that it is possible that
${\mathfrak f}[g]=\Delta_k$, for some $1 \le k \le m$).
    \item [(b)] For $j=1, \ldots, m$, we define ${\mathcal B}_j$ to be the set of
elements $[g]\in F_n/\sim$ such that ${\mathfrak f}[g]$ bounds
together with $\Delta_0$ and $\Delta_j$ an embedded pair of pants
in $M_0$.
\end{itemize}
\end{defn}

We will also need to define the half lengths, for which we need
representations into $\SLTwoC$ instead of $\PSLTwoC$. The main idea
is that by choosing a lift of the representation to $\SLTwoC$, one
can have a consistent choice of the half length for elements of
$\SLTwoC$. This follows closely the approach of Fenchel in
\cite{fenchel1989book}, and the reader is referred there for
details.

\begin{defn}\label{def:halflength}{\rm
If $\rho \in \MCSS$ and $\tilde \rho$ is a lift of $\rho$ to
$\SLTwoC$, then for an element $g \in F_n$, we define the specific
half length $l(\tilde \rho(g))/2 \in \CmodTwoPiIZ$ of $\tilde \rho
(g)$ by
\begin{eqnarray}\label{eqn:halflength}
\cosh \frac{l(\tilde \rho(g))}{2}=-\frac{{\rm tr}\,\tilde \rho(g)}{2}, %
\end{eqnarray}
with ${\mathfrak R}\,l(\tilde \rho(g))/2 >0$.}
\end{defn}

Note that the real part of the half length is just half of the
real part of the length, and both are positive, while the above
choice fixes the imaginary part, up to multiples of $2 \pi i$. The
minus sign on the right-hand side of (\ref{eqn:halflength}) is
crucial.

Our main theorem for Schottky groups can then be stated as
follows.

\begin{thm}\label{thm:mcshane schottky}
Let $\rho \in \MCSS$, and let $\tilde \rho$ be any lift of $\rho$
to $\SLTwoC$. Suppose $\rho_0$ is a fuchsian marking, with
corresponding hyperbolic surface $M_0$, and boundary components
$\Delta_0, \ldots, \Delta_m$. Let ${\mathcal P}$ and ${\mathcal
B}_j$, $j=1, \ldots, m$ be defined as in Definition
{\rm\ref{def:sets}}, relative to $M_0$. Then
\begin{eqnarray}\label{eqn:mcshane schottky}
& &\hspace{-35pt}\sum_{\{[g], [h]\}\in {\mathcal P}}
\,G\bigg(\frac{l(\tilde \rho(d_0))}{2},\frac{l(\tilde \rho(g))}{2},
\frac{l(\tilde \rho(h))}{2}\bigg)\nonumber\\
&+& \sum_{j=1}^{m}~\sum_{[g]\in {\mathcal B}_j}
S\bigg(\frac{l({\tilde \rho}(d_0))}{2},\frac{l({\tilde
\rho}(d_j))}{2},\frac{l(\tilde \rho(g))}{2}\bigg)= \frac{l({\tilde
\rho}(d_0))}{2} \mod \pi i.
\end{eqnarray}
Moreover, each series on the left-hand side of
{\rm(\ref{eqn:mcshane schottky})} converges absolutely.
\end{thm}

\begin{rmk}\label{rmk:maintheorem}\hfill {\rm
\begin{enumerate}
\item[(a)] In the case where $\rho=\rho_0$, the above is just a
reformulation of Theorem \ref{thm:complexified} for the case of a
hyperbolic surface with geodesic boundary components, and is true
without the modulo condition. In fact, the lift can be chosen so
that the right-hand side is real and positive.

\item[(b)] The identity (\ref{eqn:mcshane schottky}) is true only
modulo $\pi i$ because we have fixed the choice of the $\tanh^{-1}$
function in the definition of the functions $G(x,y,z)$ and
$S(x,y,z)$ (see Definition \ref{def:GandS}), which may differ from
the values obtained by analytic continuation by some multiple of
$\pi i$.

\item[(c)] The result is independent of the lift chosen. This is
because if $\tilde \rho$ and $\bar \rho$ are two different lifts
of $\rho$, then for each of the summands on the first series,
either ${\rm tr}\,\tilde \rho( g), {\rm tr}\,\tilde \rho( h)$ and
${\rm tr}\,\tilde \rho( d_0)$ are all equal to ${\rm tr}\,\bar
\rho( g), {\rm tr}\,\bar \rho( h)$ and ${\rm tr}\,\bar \rho( d_0)$
or exactly two of them differ by their signs (and similarly for
the summands in the second series). In the latter case, two of the
half lengths differ by $\pi i$, but it can be easily checked that
both $G(x,y,z)$ and $S(x,y,z)$ remain the same if $\pi i$ is added
to two of the arguments.

\item[(d)] The choice of the half length functions given above is
not arbitrary but arises from the computation of $G(x,y,z)$ and
$S(x,y,z)$ as ``gap'' functions (this is based on the convention
adopted by Fenchel in \cite{fenchel1989book}, see
\cite{tan-wong-zhang2004cone-surfaces} and \cite{zhang2004thesis}
for details). Roughly speaking, the relative positions of the axes
for $\tilde \rho( g)$, $\tilde \rho( h)$ and $\tilde \rho( d_0)$ are
completely determined by their traces. These axes form the
non-adjacent sides of a right angled hexagon in $\HHH$ and the half
lengths basically arise as the lengths of these sides of the
hexagon.

\end{enumerate}
}
\end{rmk}

We refer the reader to \cite{tan-wong-zhang2004schottky} for
details of the proof. We mention here that to prove the absolute
convergence of the series concerned, we use a combinatorial word
length for elements of ${\mathcal P}$ and ${\mathcal B}_j$, and by
adapting an argument from \cite{birman-series1985t}, we can show
that there is a polynomial bound (in $n$) for the number of
elements of ${\mathcal P}$ (respectively ${\mathcal B}_j$) with
combinatorial length $n$. We then show that for $\rho\in \MCSS$,
the combinatorial lengths of the elements of ${\mathcal P}$
(respectively ${\mathcal B}_j$) are comparable to the real part of
the complex lengths of their image under $\rho$ in $\PSLTwoC$ and
use these two facts to prove the absolute convergence of the
series in question.

\vskip 6pt

\begin{figure}
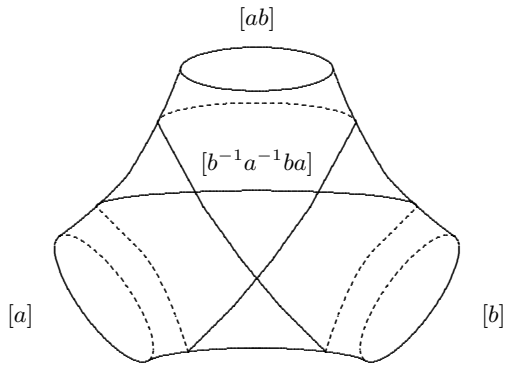

\begin{center}
\mbox{\beginpicture \setcoordinatesystem units <0.04in,0.04in> %
\setplotarea x from -40 to 40 , y from 0 to 45 %
\ellipticalarc axes ratio 3.5:1 360 degrees from 10 40  center at 0 40 %
\startrotation by .809 .588 about -20 10 \ellipticalarc axes ratio
1:3 180 degrees from -20 20 center at -20 10 \stoprotation %
\startrotation by .809 -.588 about 20 10 \ellipticalarc axes ratio
1:3 -180 degrees from 20 20 center at 20 10 \stoprotation %
\startrotation by .809 .588 about -20 10 \setdashes<1.50pt>
\ellipticalarc axes ratio 1:3 -180 degrees from -20 20 center at
-20 10 \stoprotation %
\startrotation by .809 -.588 about 20 10 \ellipticalarc axes ratio
1:3 180 degrees from 20 20 center at 20 10 \stoprotation \setsolid %
\plot -9 3 -3 9 2 15 7 22 13 33 / \plot 9 3 3 9 -2 15 -7 22 -13 33 / %
\ellipticalarc axes ratio 10:1 180 degrees from 21 22  center at 0 22 %
\setquadratic \plot -25.88 18.09 -21 22 -17 26 -13 33 -10 40 / %
\setquadratic \plot 25.88 18.09 21 22 17 26 13 33 10 40 / %
\setquadratic \plot -14.12 1.89 -9 3 0 3.5 9 3 14.12 1.89 / %
\setdashes<1.50pt> %
\ellipticalarc axes ratio 5:1 180 degrees from 13 33  center at 0 33 %
\plot -9 3 -12 11 -15 16 -18  19  -21 22 / %
\plot 9 3 12 11 15 16 18  19  21 22 / \setsolid %
\put {\mbox{\small $ [a]$}} [cb] <-1mm,-4mm> at -30 10 %
\put {\mbox{\small $ [b]$}} [cb] <3mm,-4mm> at 28 10 %
\put {\mbox{\small $ [ab]$}} [cb] <0mm,0mm> at 0 45 %
\put {\mbox{\small $ [b^{-1}a^{-1}ba]$}} [cb] <0mm,0mm> at 0 26 %
\endpicture}
\end{center}
\caption{A commutator curve on the pair of pants}\label{fig:comm}
\end{figure}

\begin{example}(A nontrivial identity for the hyperbolic pair of pants.) %
Theorem \ref{thm:mcshane schottky} can be applied to rank two
classical Schottky groups to obtain some interesting nontrivial
identities for the hyperbolic pair of pants with geodesic boundary.
The idea here is that the fundamental group in this case is free on
two generators and is isomorphic to the fundamental group of the
one-holed torus. The holonomy $\rho$ for the pair of pants is in
${\mathcal S}^{\rm mc}_{\rm alg}$, as is the holonomy  $\rho_0$ for
the one-holed torus. Using the identity obtained from $\rho_0$, one
obtains a nontrivial identity for $\rho$ via Theorem
\ref{thm:mcshane schottky}. There are interesting geometric
interpretations for each of the terms in the identity, see \S 5 of
\cite{tan-wong-zhang2004schottky} for details. Note that in this
case, the commutator $aba^{-1}b^{-1}$ of a pair of generators is a
non-simple closed curve on the pair of pants, as shown in Figure
\ref{fig:comm}, and that its trace ${\rm
tr}\,\rho(aba^{-1}b^{-1})>18$.
\end{example}

\section{{\bf The $\SLTwoC$ characters of a one-holed torus}}\label{s:character} %

The restriction of Theorem \ref{thm:complexified} to a torus with a
cusp was the original identity obtained by McShane in his thesis.
The identity (\ref{eqn:reform of cp and gb cases}) restricted to a
one-cone or one-holed torus $T$ can be regarded as generalizations
of this original identity, and reinterpreted as an identity for
representations (or more accurately, characters) of the one-holed
torus group $\pi:=\pi_1(T)$ to $\PSLTwoR$, as we saw in the previous
section. It is natural to ask how far this result can be extended to
representations into $\PSLTwoC$. Indeed one already has the
extension to quasi-fuchsian representations by Bowditch, and the
example given at the end of the previous section showed that we can
extend the identity to the case of classical Schottky
representations, including those arising from a hyperbolic pair of
pants. These, however, required some special properties including
the discreteness of the representation, which seemed unnecessarily
restrictive. For example, given a representation arising as the
holonomy of a one-cone hyperbolic torus with say cone angle $\theta$
(which is not necessarily discrete), one would expect that for
sufficiently small perturbations of the representation into
$\PSLTwoC$, the identity would still hold. This turns out to be
true, and in fact, one can give very comprehensive answers to the
questions posed above. For example, we can obtain necessary and
sufficient conditions for the generalized McShane's identity to hold
for representations/characters of $\pi$ into $\PSLTwoC$. For this,
it turns out that Bowditch's proof via algebraic/combinatiorial
methods are extremely useful, and this is the approach we use and
generalize in \cite{tan-wong-zhang2004gMm} and
\cite{tan-wong-zhang2004necsuf} to solve this problem. Another
useful corollary of this method is that we are able to obtain
various restricted and relative versions of the identity, the latter
of which have geometric interpretations in terms of punctured torus
bundles over the circle, and in particular, allows us to prove
identities for certain complete hyperbolic 3-manifolds obtained by
hyperbolic Dehn surgery on hyperbolic punctured torus bundles.

For the rest of this section, we will first start with some basic
definitions, then give statements of some of the main results, and
finally list some of the key techniques and issues involved in the
proofs of the results. We should warn the reader that the proofs are
somewhat technical in some parts, details can be found in
\cite{bowditch1998plms}, \cite{tan-wong-zhang2004gMm} and
\cite{tan-wong-zhang2004necsuf}. Note also that we shall be stating
and proving results for representations(characters) into $\SLTwoC$
instead of $\PSLTwoC$. This makes no essential difference since all
representations of $\pi$ into $\PSLTwoC$ can be lifted to $\SLTwoC$
as $\pi$ is free, and the identities obtained will be independent of
the lift chosen, and hence can be stated as identities for
$\PSLTwoC$ characters.

\vskip 6pt

\subsection{Basic Definitions}\label{ss:definitions}
Let $T$ be a one-holed torus and $\pi$ its fundamental group which
is freely generated by two elements $X,Y$ corresponding to simple
closed curves on $T$ with geometric intersection number one.

\begin{defn}\label{def:BQspace}
 The $\SLTwoC$ character variety $${\mathcal X}:={\rm
Hom}(\pi, \SLTwoC)/\!/\SLTwoC$$ of $T$ is the set of equivalence
classes of representations $\rho:\pi \mapsto \SLTwoC$, where the
equivalence classes are obtained by taking the closure of the
orbits under conjugation by $\SLTwoC$.

\end{defn}

The character variety stratifies into relative character varieties:
for $\kappa \in \CC$, the $\kappa$-relative character variety
${\mathcal X}_{\kappa}$ is the set of equivalence classes $[\rho] $
such that
$${\rm tr}\,\rho(XYX^{-1}Y^{-1})=\kappa$$ for one (and hence any)
pair of generators $X,Y$ of $\pi$. Note that the commutator
$XYX^{-1}Y^{-1}$ represents a peripheral curve in $T$.  By classical
results of Fricke, we have the following identifications:
$${\mathcal X} \cong \CC^3,$$
$${\mathcal X}_{\kappa} \cong \{(x,y,z)\in \CC^3 ~|~ x^2+y^2+z^2-xyz-2=\kappa\},$$
where the identification is given by
$$\iota:[\rho] \mapsto (x,y,z):=({\rm tr}\rho(X), {\rm tr}\rho(Y), {\rm tr}\rho(XY)),$$ %
for a fixed pair of generators $X, Y$ of $\pi$. The topology on
${\mathcal X}$ and ${\mathcal X}_{\kappa}$ will be that induced by
the above identifications.


The outer automorphism group of $\pi$, ${\rm Out}(\pi):={\rm
Aut}(\pi)/{\rm Inn}(\pi)\cong {\rm GL}(2, \mathbb Z) $ is
isomorphic to the mapping class group $\Gamma:=\pi_0({\rm
Homeo}(T))$ of $T$ and acts on ${\mathcal X}$, preserving the
trace of the commutator of a pair of generators, hence it also
acts  on ${\mathcal X}_{\kappa}$, the action is given by
$$\phi([\rho])=[\rho \circ \phi^{-1}],$$ where $\phi \in {\rm
Out}(\pi)$ and $[\rho] \in {\mathcal X}$ or ${\mathcal
X}_{\kappa}$ respectively. It is often convenient to consider only
the subgroup ${\rm Out}(\pi)^+$ of ``orientation-preserving''
automorphisms, corresponding to the orientation-preserving
homeomorphisms $\Gamma^+$ of $T$, which is isomorphic to
$\SLTwoZ$. The action of ${\rm Out}(\pi)^+$ (respectively, ${\rm
Out}(\pi)$) on ${\mathcal X}$ and ${\mathcal X}_{\kappa}$ is not
effective, the kernel is $\{\pm I\}$, generated by the elliptic
involution of $T$ so that the effective action is by $\PSLTwoZ$
(respectively, $\PGLTwoZ$).

\subsection{Simple curves; Pants graph}\label{ss:simplecurves} %

\begin{defn}\label{def:scc} We denote by
${\mathscr C}$ the set of free homotopy classes of nontrivial,
non-peripheral, unoriented simple closed curves on $T$. Elements of
${\mathscr C}$ are usually denoted by $X,Y,Z,W$.
\end{defn}
The elements of ${\mathscr C}$ correspond to certain elements of
$\pi/\! \sim$, where the equivalence relation $\sim$ is that, for
$g,h \in \pi$, $g\sim h$ if and only if $g$ is conjugate to $h$ or
$h^{-1}$. We also denote the corresponding subset of $\pi/\!\sim$ by
${\mathscr C}$, there should be no confusion.

\begin{defn}\label{def:pantsgraph}
The pants graph ${\mathscr C}(T)$  of $T$, is defined to be the
graph whose  vertices are the elements of ${\mathscr C}$, and two
vertices are joined by an edge
 if and only if the corresponding curves on $T$
have geometric intersection number one.
\end{defn}

The mapping class group $\Gamma$ and ${\rm Out}(\pi)$ act on
${\mathscr C}$ (respectively ${\mathscr C}(T)$). We can realize
${\mathscr C}(T)$ as the Farey graph/triangulation of the upper half
plane $\HH$ so that ${\mathscr C}$ is identified with $\hat
\Q:={\mathbf Q}\cup \{\infty\}$, the action of $\Gamma$ is realized
by the action of $\PGLTwoZ$ on the Farey graph. The projective
lamination space $\PL$ of $T$ is then identified with $\hat
\RR:={\mathbf R}\cup \{\infty\}$ and contains ${\mathscr C}$ as the
(dense) subset of rational points.

\vskip 5pt

\subsection{Bowditch Q-conditions (BQ-conditions)}\label{ss:BQconditions}
 We define a certain subspace of ${\mathcal X}$ which we will call the Bowditch space.
 First note that for $[\rho] \in {\mathcal X}$ and $X \in {\mathscr C}$,
${\rm tr}\, \rho(X)$ is well-defined.

\begin{defn} \label{def:BQspace}

The \emph{Bowditch space} is   the subset ${\mathcal
X}_{BQ}\subset {\mathcal X}$ consisting of characters $[\rho]$
satisfying the following conditions (the \emph{Bowditch
Q-conditions}):
\begin{enumerate}
    \item ${\rm tr}\,\rho(X) \not \in [-2,2]$ for all $X \in \mathscr C$;
    \item $|{\rm tr}\,\rho(X)| \le 2$ for only finitely many (possibly no)
    $X \in \mathscr C$.
\end{enumerate}
\end{defn}

 For a fixed $[\rho] \in  {\mathcal X}$ and $U \subset \mathscr
C$, we say that the BQ-conditions are satisfied on $U$ for
$[\rho]$ if conditions (1) and (2) above hold for all $X \in U$.

\subsection{Statement of results for $\SLTwoC$ characters}\label{ss:resultscharacters}
\vskip 5pt We have the following extension and generalization of
Theorem \ref{thm:complexified} to  characters in ${\mathcal X}$.

\begin{thm}\label{thm:TWZ}{\rm(Theorems 2.2, 2.3 and Proposition 2.4 of
\cite{tan-wong-zhang2004gMm})}

    {\rm(a)} Bowditch space ${\mathcal X}_{BQ}$ is open in the whole character space ${\mathcal X}$.

    {\rm(b)} The mapping class group $\Gamma$ acts properly discontinuously on ${\mathcal X}_{BQ}$.
Furthermore, ${\mathcal X}_{BQ}$ is the largest open subset of
${\mathcal X}$ for which this holds.

    {\rm(c)} For a character $[\rho] \in {\mathcal X}_{BQ} \cap {\mathcal X}_{\kappa}$,
\begin{eqnarray}\label{eqn:TWZ}
\sum_{X \in {\mathscr C}}\log\frac{e^{\nu}+e^{l(\rho(X))}}{e^{-\nu}+e^{l(\rho(X))}} =\nu \mod 2 \pi i, %
\end{eqnarray}
where $\nu=\cosh^{-1}(-\kappa/2)$, and the sum converges absolutely.
\end{thm}

\vskip 6pt

\begin{rmk}\label{rmk:1}~\\
\vspace{-12pt}
\begin{enumerate}
    \item The (complex) length $l(\rho(X))$ is related to the trace as in equation
    (\ref{eqn:complexlength}).
    \item We are using the formula for $G(x,y,z)$ given in
    (\ref{eqn:G in log}) for part (c), note that there are no $S(x,y,z)$ terms
    since there is only one boundary component.
    \item In the case when $\kappa=-2$, $\nu=0$ and all the terms of
(\ref{eqn:TWZ}) are identically zero. However, if we take the first
order infinitesimals, or the formal derivative of (\ref{eqn:TWZ})
with respect to $\nu$ and evaluate at $\nu=0$, we get
\begin{equation}\label{eqn:McShanetorus}
\sum_{X \in {\mathscr C}}\frac{1}{1+e^{l(\rho(X))}} =\frac 12,
\end{equation}
which is McShane's original identity in \cite{mcshane1991thesis}
for real type-preserving characters, and also Bowditch's
generalization in \cite{bowditch1996blms} and
\cite{bowditch1998plms} for  type-preserving characters satisfying
the BQ-conditions.
    \item When $\kappa=2$, which corresponds to the reducible
characters,  the identity is also trivial. In this case, however,
the Bowditch Q-conditions are never satisfied, see
\cite{tan-wong-zhang2004endinvariants}.
    \item Parts (a) and (b) of the above were originally
stated in \cite{tan-wong-zhang2004gMm} in terms of the relative
character varieties ${\mathcal X}_{\kappa}$.
    \item $\nu$ is a specific choice of half of the complex length
    of the peripheral curve on $T$, note that the minus sign is
    crucial for the identity to hold.
\end{enumerate}
\end{rmk}

\subsection{Necessary and sufficient conditions}\label{ss:necandsuff}

    Replacing condition (1) of the BQ-conditions by ($1'$) ${\rm
tr}\,\rho(X) \not \in (-2,2)$ for all $X \in \mathscr C$, we get the
\emph{extended Bowditch space} $\hat {\mathcal X}_{BQ}$, and we have
the following result:
\begin{thm} \label{thm:necsuff} {\rm(Theorem 1.5 of
\cite{tan-wong-zhang2004necsuf})} For $[\rho] \in {\mathcal X}$, the
identity {\rm(\ref{eqn:TWZ})} of Theorem {\rm\ref{thm:TWZ}(c)} holds
{\rm(}with absolute convergence of the sum{\rm)} if and only if
$[\rho]$ lies in the extended Bowditch space $\hat {\mathcal
X}_{BQ}$.
\end{thm}

The above result gives a complete answer to the question of when
the generalized McShane's identity holds for $\SLTwoC$ characters
of $T$.

\subsection{McShane-Bowditch identities for punctured torus bundles}%
\label{ss:McShaneBowdithc}

We next consider further variations of the McShane-Bowditch
identities. Recall that $\theta \in {\rm Out}(\pi) \cong\Gamma$ acts
on ${\mathcal X}$ where the action is given by
$$\theta([\rho])=[\rho \circ \theta^{-1}].$$

Suppose that $[\rho] \in {\mathcal X}$ is stabilized by an Anosov
element $\theta \in \Gamma^+$ (this corresponds to a hyperbolic
element if we identify $\Gamma^+$ with $\SLTwoZ$), that is,
$\theta([\rho])=[\rho]$. We can associate to this a representation
of $\pi_1(M$) into $\SLTwoC$, where $M$ is a punctured torus bundle
over the circle, with monodromy $\theta$. The restriction of the
representation to the fibre is $[\rho]$. We can find a specific lift
of $\theta$ to ${\rm Aut}(\pi)$ which corresponds to choosing a
specific longitude of the boundary torus of $M$ (see
\cite{bowditch1997t} or \cite{tan-wong-zhang2004gMm} for details).
So fixing a representation $\rho$ in the class $[\rho]$, there
exists $A \in \SLTwoC$ such that for all $\alpha \in \pi$,
$$\theta(\rho)(\alpha)=A\cdot \rho(\alpha)\cdot A^{-1}.$$
Note that ${\rm tr}\,A$ is independent of the choice of $\rho$ in
the conjugacy class $[\rho]$. Note also that ${\rm tr}\,\rho(X)$
is well-defined on the equivalence classes $[X] \in \mathscr
C/\langle \theta \rangle $. Suppose further that $[\rho]$
satisfies the \emph{relative} Bowditch Q-conditions on ${\mathscr
C}/\langle \theta \rangle$, that is,
\begin{enumerate}
    \item ${\rm tr}\,\rho(X) \not \in [-2,2]$ for all %
    $[X] \in \mathscr C/\langle \theta \rangle $;
    \item $|{\rm tr}\,\rho(X)| \le 2$ for only finitely many
    $[X] \in \mathscr C/\langle \theta \rangle $.
\end{enumerate}
Using the identification of $\Gamma^+$ with $\SLTwoZ$ and
${\mathscr C}$ with $\hat \Q \subset \hat \RR \cong {\mathscr
{PL}}$ in \S \ref{ss:simplecurves},  we get that the repelling and
attracting fixed points of $\theta$, $\mu_-, \mu_+ \in {\mathscr
{PL}}$ partition ${\mathscr C}$ into two subsets ${\mathscr C}_L
\sqcup {\mathscr C}_R$ which are invariant under the action of
$\theta$. We have the following generalizations of the
McShane-Bowditch identities:
\begin{thm}\label{thm:bundle} {\rm (Theorems 5.6 and 5.9 of
\cite{tan-wong-zhang2004gMm})} Suppose that $[\rho]$ is stabilized
by an Anosov element $\theta \in \Gamma^+$ and satisfies the
relative Bowditch Q-conditions as stated above. Then
\begin{eqnarray}\label{eqn:bundle1}
\hspace{-20pt}\sum_{[X] \in {\mathscr C}/\langle \theta
\rangle}\log\frac{e^{\nu}+e^{l(\rho(X))}}{e^{-\nu}+e^{l(\rho(X))}}
=0 \mod 2 \pi i,
\end{eqnarray}
and
\begin{eqnarray}\label{eqn:bundle2}
\sum_{[X] \in {\mathscr C}_L/\langle \theta
\rangle}\log\frac{e^{\nu}+e^{l(\rho(X))}}{e^{-\nu}+e^{l(\rho(X))}}
=\pm l(A) \mod 2 \pi i,
\end{eqnarray}
where the sums converge absolutely; and $l(A)$ is the complex length
of the conjugating element $A$ corresponding to $\theta$ as
described above, and the sign in {\rm(\ref{eqn:bundle2})} depends
only on our choice of orientations.
\end{thm}

\begin{rmk}
For type-preserving characters ($\kappa=-2$), the result is due to
Bowditch \cite{bowditch1997t}, where the summands of
(\ref{eqn:bundle1}) and (\ref{eqn:bundle2}) should be replaced
appropriately as in Remark \ref{rmk:1}(3) by the summands of
McShane's original identity, and $l(A)$ in (\ref{eqn:bundle2})
should be replaced by $\lambda$,  the modulus of the cusp of $M$
with the complete, finite volume hyperbolic structure. There are
also similar identities in the case where $\theta$ is reducible,
that is, corresponds to a parabolic element of $\SLTwoZ$, see
\cite{tan-wong-zhang2004necsuf}.
\end{rmk}

The above result has applications to closed hyperbolic 3-manifolds.
As before, let $M$ be an orientable 3-manifold which fibers over the
circle, with the fiber a once-punctured torus, $T$ and suppose that
the monodromy $\theta$ of $M$ is Anosov. By results of Thurston, see
\cite{thurston19??am} and \cite{thurston1978notes}, $M$ has a
complete finite-volume hyperbolic structure with a single cusp,
which can in turn be deformed to incomplete hyperbolic structures,
on which hyperbolic Dehn surgery can be performed to obtain complete
hyperbolic manifolds without cusps. Restricting the holonomy
representation to the fiber gives us characters which are stabilized
by $\theta$, and in the complete case, the relative Bowditch
Q-conditions are satisfied (see \cite{bowditch1997t}). For small
deformations of the complete structure to incomplete structures, the
relative BQ-conditions are still satisfied since these are open
conditions (see \cite{tan-wong-zhang2004gMm}). The identities can be
interpreted as series identities for these (in)complete structures,
involving the complex lengths of certain geodesics corresponding to
the homotopy classes of essential simple closed curves on the fiber.
The quantity $\nu$ can be interpreted as half the complex length of
the meridian of the boundary torus, and $l(A)$ as the complex length
of a (suitably chosen) longitude of the boundary torus. In
particular, the identity can be interpreted as an identity for the
closed hyperbolic 3-manifolds obtained by hyperbolic Dehn surgery on
the original complete manifold, if the Dehn surgery invariants are
sufficiently close to $\infty$.


One question which arises is whether the identity holds for all
closed hyperbolic 3-manifolds obtained by hyperbolic Dehn surgery on
a hyperbolic punctured torus bundle over the circle. The openness of
the relative BQ-conditions ensures that this is true for almost all
such manifolds (except for possibly a finite number of exceptions).
Another question arising is whether a similar result holds in the
case of hyperbolic Dehn surgery on punctured surface bundles, as
studied by Akiyoshi, Miyachi and Sakuma
\cite{akiyoshi-miyachi-sakuma2004preprint}.

\subsection{Key points used in the proofs }\label{ss:keypoints}
The combinatorial structure of $\C(T)$ as well as the Fricke trace
relation which can be interpreted as an edge relation play
fundamental roles which we sketch here.

Recall that $\C(T)$ has the structure of the Farey tessellation, and
the set of vertices $\C$ can be identified with $\hat \Q$. The dual
graph $\Sigma$ to $\C(T)$ is a trivalent tree whose complementary
regions can be identified with the vertices of $\C(T)$. Denote by
$V(\Sigma)$, $E(\Sigma)$, $\vec E(\Sigma)$ and $\Omega(\Sigma)$ the
sets of vertices, edges, directed edges and complementary regions of
$\Sigma$ respectively. Call $(X,Y) \in \C \times \C$ a generating
pair if $X$ and $Y$ are connected by an edge in $\C(T)$, and
$(X,Y,Z)\in \C \times \C \times \C$ a generating triple if $X,Y$ and
$Z$ are the vertices of a triangle in $\C(T)$. Generating pairs
correspond to edges of $\Sigma$ and generating triples correspond to
vertices of $\Sigma$. More specifically, to an edge $e$ of $\Sigma$,
we write $e = (X,Y;Z,Z')$ if $(X,Y)$ corresponds to $e$ and
$(X,Y,Z)$, $(X,Y,Z')$ are generating triples. Similarly, we use
$\vec e = (X,Y;Z \rightarrow Z')$ to indicate that the directed edge
$\vec e$ points from $Z$ to $Z'$, see Figure \ref{fig:edge
oriented}, where we have drawn part of $\Sigma$, and used the
identification of $\Omega(\Sigma)$ with $\C$. Denote by $-\vec e$
the directed edge with the opposite direction to $\vec e$. For $\vec
e = (X,Y;Z \rightarrow Z')$, we define $\textrm{Tail}(\vec e)$, the
tail of $\vec e$ to be the subset of $\C$ in the interval between
$X$ and $Y$ (inclusive) which contains $Z$. In particular,
$\textrm{Tail}(\vec e)\cup \textrm{Tail}(-\vec e)=\C$, and
$\textrm{Tail}(\vec e) \cap \textrm{Tail}(-\vec e)=\{X,Y\}$.

For each character $[\rho] \in {\mathcal X}_{\kappa}$, by taking the
trace function, we obtain a {\em trace map} 
\begin{eqnarray*}
\phi: \C \rightarrow \CC \quad \textrm{where} \quad \phi(X)={\rm tr}\rho(X). %
\end{eqnarray*}
(We call it a generalized Markoff map in
\cite{tan-wong-zhang2004gMm} following \cite{bowditch1998plms}.)

Henceforth, for a fixed trace map $\phi$, we adopt the convention of
using the lower case letters to represent the values of $\phi$, that
is, $\phi(X)=x$, $\phi(Y)=y$, etc. Then $\phi$ satisfies the
following vertex and edge relations, arising from the Fricke trace
identities: \vskip 3pt

{\bf Vertex relation}. For every generating triple $(X,Y,Z)$, %
\begin{equation}\label{eqn:vertex}
    x^2+y^2+z^2-xyz-\kappa-2=0.
\end{equation}

{\bf Edge relation}. For every edge $e=(X,Y;Z,Z')$, %
\begin{equation}\label{eqn:edge}
    z+z'=xy.
\end{equation}

It turns out that the edge relation is more fundamental than the
vertex relation. To start with, one can show easily that if the
edge relation is satisfied for all edges, than the vertex relation
propagates along the edges to cover the entire tree $\Sigma$.
Secondly, $\phi$  is completely determined by its values on any
generating triple $(X,Y,Z)$ by successively applying the edge
relation (\ref{eqn:edge}).

Each $[\rho] \in {\mathcal X}$ (equivalently, the induced trace
map $\phi$ on $\C$) determines a map $f:E(\Sigma) \rightarrow \vec
E(\Sigma)$, where each edge $e$ is assigned a direction or flow
from the larger absolute value to the smaller one, that is,
$$f(e)=\vec e = (X,Y;Z \rightarrow Z')$$
if $|z|\ge |z'|$. There is some ambiguity when $|z|=|z'|$ in which
case we can assign either direction. This ambiguity does not affect
the large scale behavior of $f(E(\Sigma))$, except in some very
special trivial cases. We then have the following elementary but
important results (see \cite{tan-wong-zhang2004gMm} for proofs). For
the purposes of our discussion, we fix $[\rho]\in {\mathcal
X}_{\kappa}$ where $\kappa \neq 2$ ($[\rho]$ is not reducible), with
corresponding trace map $\phi$.

\begin{prop}\label{prop:twoout}
If $(X,Y,Z)$ is a generating triple corresponding to the vertex $v
\in V(\Sigma)$ and $f(e)$ points away from $v$ for at least two of
the edges adjacent to $v$, then ${\rm min}(|x|,|y|,|z|) \le 2$.
\end{prop}
%
\begin{lem}\label{lem:fundamental}{\rm(Bowditch \cite{bowditch1998plms})} %
For all $K \ge 2$,  $\C(K):=\{X\in \C \mid \phi(X) \le K\}$ is
connected, that is, the subgraph of $\C(T)$ spanned by $\C(K)$ is
connected. In particular, $\C(2)$ is connected.

\end{lem}

The above can be regarded as a quasi-convexity result, namely, for
any $K \ge 2$, for any $X, Y \in \C(K)$, the geodesic in $\C(T)$
joining $X$ to $Y$ is a bounded distance from the subgraph in
$\C(T)$ spanned by $\C(K)$.

\begin{prop}\label{prop:neighbors}
Suppose that $X \in \C$ and $Y_n$, $n \in {\mathbf Z}$ are the
neighbors of $X$, in cyclical order.

    {\rm(a)} If $x \not\in [-2,2] \cup \{\pm \sqrt{\kappa+2}\}$, then %
$\lim_{n\rightarrow \pm \infty}|y_n|=\infty$ with exponential growth in $|n|$. %

    {\rm(b)} If $x=\pm 2$ and $\kappa \neq 2$, then %
$\lim _{n\rightarrow \pm\infty}|y_n|=\infty$ with linear growth in $|n|$. %
\end{prop}

The proof of Theorem \ref{thm:TWZ} now proceeds as follows:

First we show that the BQ-conditions are open conditions. This is
achieved by showing that the conditions are controlled by a finite
subtree of $\Sigma$; the proof is essentially that given by Bowditch
in \cite{bowditch1998plms}, with some slight modifications.

\begin{figure}
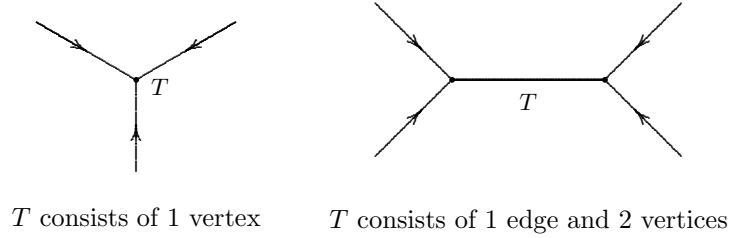

\begin{center}
\mbox{\beginpicture \setcoordinatesystem units <0.6in,0.6in> %
\setplotarea x from -1.2 to 1.2, y from -1.2 to 1.2 %
\plot 0 -0.8 0 0 -0.866 0.5 / %
\plot 0 0 0.866 0.5 / %
\arrow <6pt> [.16,.6] from 0 -0.5 to 0 -0.4 %
\arrow <6pt> [.16,.6] from -0.866 0.5  to -0.433 0.25 %
\arrow <6pt> [.16,.6] from 0.866 0.5  to 0.433 0.25 %
\put {\mbox{\Huge $\cdot$}} [cc] <0mm,-0.2mm> at 0 0 %
\put {\mbox{$T$ consists of 1 vertex}} [ct] <0mm,0mm> at 0 -1.12 %
\put {\mbox{\small $T$}} [lc] <2mm,-1mm> at 0 0 %
\endpicture} %
\hspace{0.2in} \raisebox{-11ex} %

\mbox{\beginpicture \setcoordinatesystem units <0.4in,0.4in> %
\setplotarea x from -1.2 to 3.2, y from -1.2 to 1.2 %
\plot -1 1 0 0 2 0 3 1 / %
\plot -1 -1 0 0 / \plot 0 0.01 2 0.01 / %
\plot 0 -0.01 2 -0.01 / \plot 2 0 3 -1 / %
\arrow <6pt> [.16,.6] from  -0.5 0.5  to -0.4 0.4 %
\arrow <6pt> [.16,.6] from -0.5 -0.5 to -0.4 -0.4 %
\arrow <6pt> [.16,.6] from 2.5 0.5  to 2.4 0.4 %
\arrow <6pt> [.16,.6] from 2.5 -0.5  to 2.4 -0.4 %
\put {\mbox{\Huge $\cdot$}} [cc] <0mm,-0.2mm> at 0 0 %
\put {\mbox{\Huge $\cdot$}} [cc] <0mm,-0.2mm> at 2 0 %
\put {\mbox{$T$ consists of 1 edge and 2 vertices}} [ct] <0mm,0mm> at 1 -1.7 %
\put {\mbox{\small $T$}} [cc] <0mm,-3mm> at 1 0 %
\endpicture} %
\end{center} %
\caption{Two simple circular sets $C(T)$}\label{fig:circular} %
\end{figure}

Next, we show that if $\phi$ satisfies the BQ-conditions, then the
function $\log^+|\phi|:=\max\{0,\log|\phi|\}$ on $\C$ has lower
Fibonacci growth. Roughly, this means that it is comparable to the
combinatorial (word) length function on $\C$, that is, there exists
some $k>0$ such that $\log^+|\phi(X)| \ge k\|X\|_w$ for all but a
finite number of $X \in \C$, where $\|X\|_w$ is the (cyclically
reduced) word length of $X$ with respect to any fixed  pair  of
generators for $\pi$. This is achieved by showing that one can find
a finite subtree $T$ of $\Sigma$ such that $f(e)$ is directed
towards $T$ for all $e \in \Sigma \setminus T$, by applying
Proposition \ref{prop:twoout}, Lemma \ref{lem:fundamental} and
Proposition \ref{prop:neighbors}. We can think of this subtree $T$
as the union of sufficiently long boundary paths of the elements of
$\C(2)$ in $\Sigma$ (which by the BQ-conditions is finite). Define
the circular set $C(T)$ of $T$ to be the set of directed edges $\vec
e \in \vec E(\Sigma)$ adjacent to $T$ and directed towards $T$, see
Figure \ref{fig:circular} where examples of $C(T)$ are given for the
two simplest cases of $T$. Then $C(T)$ is finite, and $\bigcup_{\vec
e \in C(T)} \textrm{Tail}(\vec e)=\C$. The lower bound is achieved
by showing that the lower bound holds for each of the sets
$\textrm{Tail}(\vec e)$ (with possibly different constants $k$),
where $\vec e \in C(T)$, from which the general result on $\C$
follows since it is a finite union of these sets. Now Theorem
\ref{thm:TWZ}(b) follows from the lower Fibonacci growth, and one
deduces in a fairly straightforward manner that the series in
(\ref{eqn:TWZ}) of Theorem \ref{thm:TWZ}(c) converges absolutely,
recalling the relation between the traces and the complex lengths
given in (\ref{eqn:complexlength}).

To complete the proof of Theorem \ref{thm:TWZ}(c) and show that the
sum is indeed as stated, we need to modify the methods of Bowditch
\cite{bowditch1998plms} slightly. First note that $\kappa=-2$ for
all the characters $[\rho]$ considered in \cite{bowditch1998plms}
(also called type-preserving). Bowditch used the edge-weight
function (which depends on the character $[\rho]$, or, equivalently,
the corresponding trace map $\phi$) %
\, $\psi:=\psi_{\phi} : \vec E(\Sigma) \rightarrow {\mathbf C}$ %
defined by
\begin{equation}\label{eqn:xy/z}
    \psi(\vec e)=z/xy,
\end{equation}
where \,$\vec e=(X,Y;Z'\rightarrow Z)$, and $x=\phi(X),
~y=\phi(Y),~z=\phi(Z)$ in a very ingenious manner to prove the
result, by taking sums of $\psi(\vec e)$ over circular sets of
larger and larger subtrees $T$ of $\Sigma$ which eventually exhaust
$\Sigma$. The main properties of $\psi$ which were used are the
following:

\begin{lem}\label{lem:psi} {\rm(Properties of the edge-weight function
$\psi$ in the case $\kappa=-2$.)}\,\, Suppose $\phi$ corresponds to
a type-preserving character $[\rho]$ {\rm(}that is,
$\kappa=-2${\rm)} and $\phi(X) \neq 0$ for all $X \in \Omega$.
If $\psi:=\psi_{\phi}$ is the edge-weight function defined 
by {\rm(\ref{eqn:xy/z})}, then

{\rm (i)} for a directed edge $\vec e \in \vec E$,
\begin{eqnarray}\label{eqn:sum psi i}
\psi(\vec e)+\psi(-\vec e)=1;
\end{eqnarray}

{\rm (ii)} for a circular set $C(T) \subset \vec E$,
\begin{eqnarray}\label{eqn:sum psi ii}
\textstyle\sum_{\vec e \in C(T)} \psi(\vec e)=1.
\end{eqnarray}
\end{lem}

Note that (\ref{eqn:sum psi i}) is just the edge relation
(\ref{eqn:edge}), and (\ref{eqn:sum psi ii}) is just the vertex
relation (\ref{eqn:vertex}) in the case where $\kappa=-2$ and $T$ a
vertex. The case for general $T$ follows easily by an inductive
argument using (\ref{eqn:sum psi i}).


 Now for general $[\rho]\in {\mathcal X}_{\kappa}$ where $\kappa \neq
\pm 2$, we need to find a corresponding edge weight function $\psi$
(depending on $[\rho]$ and $\kappa$) such that a suitable
generalization of Lemma \ref{lem:psi} holds. It turns out we can
define an edge-weight function $\psi$ such that:

\begin{lem}\label{lem:psi2} {\rm(Properties of the general edge-weight function
$\psi$ for $\kappa \neq \pm 2$)}\,\, Suppose $[\rho]\in {\mathcal
X}_{\kappa}$, $\kappa \neq \pm 2$,  with corresponding trace map
$\phi$, and $\phi(X) \neq 0, \pm \sqrt{\kappa+2}$ for all $X \in
\C$. Then

{\rm (i)} for a directed edge $\vec e \in \vec E$,
\begin{eqnarray}\label{eqn:sum psi i_2}
\psi(\vec e)+\psi(-\vec e)=\nu \mod 2 \pi i;
\end{eqnarray}

{\rm (ii)} for a circular set $C(T) \subset \vec E$,
\begin{eqnarray}\label{eqn:sum psi ii_2}
\textstyle\sum_{\vec e \in C(T)} \psi(\vec e)=\nu \mod 2 \pi i,
\end{eqnarray}
\noindent where $\nu=\cosh^{-1}(-\kappa/2)$ and for $\vec
e=(X,Y;Z'\rightarrow Z)$
\begin{eqnarray}\label{psi(x,y,z)=log}
\psi(\vec e):= \log\bigg(\frac{1+(e^{\nu}-1)(z/xy)}%
{\sqrt{1-(\kappa+2)/x^2}\sqrt{1-(\kappa+2)/y^2}}\bigg). %
\end{eqnarray}
\end{lem}

The rest of the proof of Theorem \ref{thm:TWZ}(c) then proceeds
along the same lines as that used in \cite{bowditch1998plms}, with
some extra technicalities required to prove that the error term
approaches zero. It should be pointed out that the function $\psi$
has a geometric interpretation in the case of a hyperbolic one-cone
torus in terms of certain gaps formed by the simple geodesics
emanating and terminating  at the cone point (c.f. \S \ref{s:intro}). %
We have a similar geometric interpretation for $\psi$ for a
hyperbolic one-holed torus with geodesic boundary. In these cases,
the properties (\ref{eqn:sum psi i_2}) and (\ref{eqn:sum psi ii_2})
are almost self-evident from the geometric interpretation. The
general formula for $\psi$ is then obtained from these geometric
cases by analytic continuation.


For Theorem \ref{ss:McShaneBowdithc}, one requires  restricted and
relative versions of Theorem \ref{thm:TWZ}. These can be obtained
fairly easily now that we have found the edge-weight function
$\psi$ with the desired properties. The proofs follow essentially
the same lines as Bowditch's in \cite{bowditch1997t}. Details can
be found in \cite{tan-wong-zhang2004gMm} or
\cite{tan-wong-zhang2004necsuf}.


A challenge is to extend the methods described above to give similar
algebraic/combinatorial proofs for characters of general punctured
surfaces. It seems that a good understanding of the pants graph or
curve complex for such surfaces is essential, although this was
circumvented rather cleverly in the work of Akiyoshi, Miyachi and
Sakuma in \cite{akiyoshi-miyachi-sakuma2004preprint} when they
studied punctured torus bundles over the circle.


For characters (of the one-holed torus) which do not satisfy the
(extended) BQ-conditions, or the relative BQ-conditions, the
dynamics of the action of the mapping class group is very
interesting, and the general methods described above can be applied
to study the question. For a character $[\rho]\in \X$, we say that
$\lambda \in \PL$ is an {\it end invariant} of $[\rho]$ if there
exists a sequence of distinct $X_k \in \C$ converging to $\lambda$
such that $|{\rm tr}\rho(X_k)|$ is bounded. It is easy to see that
the set of end invariants is empty if $[\rho]\in \hat {\mathcal
X}_{BQ}$ by Theorem \ref{thm:necsuff}, and is equal $\{\mu_+,
\mu_-\}$ if $[\rho]$ satisfies the conditions of Theorem
\ref{thm:bundle}. In particular, the conjecture is that the set of
end invariants is a Cantor set if it contains at least three
elements and is not the entire projective lamination space $\PL$.
See \cite{tan-wong-zhang2004endinvariants}, where the set of end
invariants was studied in various cases, with supporting evidence
for the conjecture.

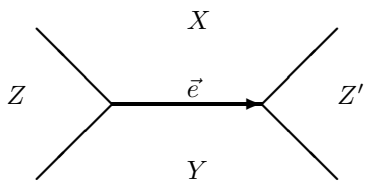
\begin{figure}
\setlength{\unitlength}{1mm} 
\begin{picture}(60,30)
\thicklines \put(20,10){\vector(1,0){20}}
\put(20,10){\line(-1,1){10}} \put(20,10){\line(-1,-1){10}}
\put(40,10){\line(1,1){10}} \put(40,10){\line(1,-1){10}}
\put(30,20){$X$} \put(30,0){$Y$} \put(50,10){$Z'$} \put(6,10){$Z$}
\put(30,11){$\vec e$}
\end{picture}
\caption{The directed edge $\vec e = (X,Y;Z \rightarrow Z')$}
\label{fig:edge oriented}
\end{figure}



\end{document}